\documentclass[11pt]{article}

\usepackage{amsmath, amsthm, amssymb}
\usepackage{setspace}
\usepackage{graphicx}
\usepackage{xcolor}
\newtheorem {definition}{Definition}
\newtheorem {remark}{Remark}

\newtheorem {theorem}{Theorem}
\newtheorem*{thm1}{Theorem 1}
\newtheorem*{thm2}{Theorem 2}
\newtheorem {lemma}{Lemma}

\newtheorem {proposition}{Proposition}

\numberwithin{equation}{section}

\title{\textbf{A model in one-dimensional thermoelasticity}}

\author{
	Tomasz Cie\'slak\thanks{Institute of Mathematics, Polish Academy of Sciences, Poland, \texttt{cieslak@impan.pl}},  
	Marija Gali\' c\thanks{Department of Mathematics, Faculty of Science, University of Zagreb, Croatia, \texttt{marija.galic@math.hr}}, 
	Boris Muha\thanks{Department of Mathematics, Faculty of Science, University of Zagreb, Croatia, \texttt{borism@math.hr}}
	}
\date{}

\begin{document}

\maketitle
\begin{abstract}
We study a one-dimensional nonlinear hyperbolic-parabolic initial boundary 
value problem occurring in the theory of thermoelasticity.
We prove existence and uniqueness of the local-in-time strong solution. 
Also, some global-in-time weak measure valued solutions are proven to exist. 
To this end we introduce an auxiliary 
problem with artificial viscosity and prove its global-in-time well-posedness. Next, 
we show that solutions of the auxiliary problem converge, at some short time interval 
to the strong solution, and to our measure valued solution
for an arbitrary time.
\end{abstract}

\section{Introduction}\label{sec:intro}
In the thermoelasticity theory, it is assumed that a considered material is elastic 
and its response to an external load depends on the temperature.
We consider the following system which is the simplest model of  
the one-dimensional thermoelasticity which takes into account the 
nonlinear coupling between the temperature and displacement of the 
elastic material (see \cite{DafHsiao,Daf,rashi,slemrod} for related 1D problems):
\begin{equation}\label{Evolution}
\begin{array}{c}
u_{tt}-u_{xx}=-\mu \theta_x \quad {\rm in}\quad (0,\infty)\times (a,b), \\
\theta_t-\theta_{xx}+\mu\theta(u_t)_x=0,\; \theta\geq 0
\quad{\rm in}\quad (0,\infty)\times (a,b),\\
u(t,a)=u(t,b)=0,\; \theta_x(t,a)=\theta_x(t,b)=0,\\
u(0,.)=u_0,\; u_t(0,.)=u_1,\\
\theta(0,.)=\theta_0\geq 0.
\end{array}
\end{equation}
In this system, $u$ denotes the displacement of the elastic material, $\theta$
is the material's temperature and $\mu$ is a material's constant.
We are interested in positive solutions $\theta\geq 0$. 
The coupling term $\mu\theta(u_t)_x$ in the temperature equation enables us to
use the comparison principle at a formal level and arrive at non-negativity of the temperature, 
unlike in the linear approximation where the coupling term is replaced
by $\mu (u_t)_x.$

Let us briefly describe the derivation of the system \eqref{Evolution} from 
the first principles.
The balance of momentum for an elastic material reads
\begin{equation}\label{wave1}
\rho u_{tt} = \text{div}\sigma + f
\end{equation}
where $u$ is the displacement of the material, $\rho$ is the mass density, $\sigma$ is the Piola-Kirchhoff
stress tensor and $f$ describes the density of external  forces acting on the continuum. 
The linear elastic constitutive stress-strain relation is a generalization
of the classical Hooke's law
\begin{equation}\label{wave2}
\sigma = \mathcal{D}\text{sym}(\nabla u)-\mu(\theta-\theta_*)I,
\end{equation}
where $\mathcal{D}$ is the elasticity tensor (which is symmetric and positive definite),
$\text{sym}(\nabla u) = \frac{1}{2}(\nabla u + \nabla^T u)$ is the symmetrized gradient of $u$
and $\mu$ is a positive constant which depends on the material and is determined experimentally.
The material's temperature is denoted by $\theta,$ 
while $\theta_*$ denotes a given, reference temperature. 
Assuming that there is no external load acting to the material, and
by taking the elasticity tensor $\mathcal{D}$ to be equal to the identity and
the mass density $\rho$ to be equal to 1, the system \eqref{wave1}-\eqref{wave2}
in 1D reduces exactly to \eqref{Evolution}$_1.$

In the thermoelasticity theory, the system of elastic equations \eqref{wave1}-\eqref{wave2} is coupled
with the heat equation which describes the evolution of the material's temperature,
and is a consequence of \emph{the first law of thermodynamics}:  the time derivative
of the total energy is equal to the sum of the power of external forces and the
rate of heat received by the continuum
$$
\frac{dE}{dt}(t)=P_{ext}(t)+\frac{dQ}{dt}(t).
$$
This principle implies the differential equation
$$
\frac{de}{dt}(t) +\text{div} q=\sigma\cdot\nabla v,
 $$
where $e$ is the density of the internal energy, $q$ is the heat flux and $v$ is the velocity of the material's displacement.
Assuming that $e=c\theta + \mathcal {D}^{-1}T\cdot T,$  where $c$ is a constant
depending on the material and $T=\mathcal{D}\text{sym}(\nabla u),$ and, furthermore,
assuming that the heat transfer 
satisfies the Fourier law $q=-\kappa\nabla\theta,$ with $\kappa$ being the material's
conductivity, we obtain the heat equation in the form 
\begin{equation}\label{heat}
c\theta_t - \kappa \Delta \theta + \mu (\theta - \theta_*)
\text{div} u_t = 0.
\end{equation}
Assuming all the physical constants to be equal to 1, i.e. $c=\kappa=1,$ 
the reference temperature $\theta_*$ to be equal to zero, the equation \eqref{heat} in 1D reduces to \eqref{Evolution}$_1.$ 

We prove the existence and uniqueness, locally in time,
of a strong solution, with the temperature being \emph{non-negative},
to this hyperbolic-parabolic coupling problem by using 
the artificial viscosity approach. We would like to emphasize that, 
physically, temperature is expected to be non-negative but the models
studied so far do not seem to fulfill this expectation.
 
To be more precise, we first introduce an auxiliary,
regularized problem, and passing to the limit as regularization parameter tends to zero, 
we obtain a solution to the original problem \eqref{Evolution}.
Furthermore, we also prove the existence, globally in time, of a measure valued solution
to problem \eqref{Evolution}, as defined below. These results are summarized in the following two
theorems:

\begin{theorem}\label{thm:exist_local}
	Let $u_0\in H^2(a,b)\cap H^1_0(a,b)$, $u_1\in H_0^1(a,b)$, $\theta_0\in H^2(a,b)$, $\theta_0\geq 0$. 
	Then there exists a time $T_0>0$ (depending on data) and a unique solution $(u,\theta)$ 
	to problem \eqref{Evolution} on $(0,T_0)$ with the following regularity:
	\begin{align*}
	&\|u\|_{W^{2,\infty}(0,T_0;L^2(a,b))}+\|u\|_{W^{1,\infty}(0,T_0;H^1(a,b))}
	+\|u\|_{L^{\infty}(0,T_0;H^2(a,b))}\\
	&\quad+\|\theta\|_{W^{1,\infty}(0,T_0;L^2(a,b))}+\|\theta\|_{H^1(0,T_0;H^1(a,b))}\leq C.
	\end{align*}
\end{theorem}

\iffalse
\begin{theorem}\label{thm:exist_global}
	Let $u_0\in H_0^1(a,b)$, $u_1\in L^2(a,b)$, $\theta_0\in L^2(a,b)$, $\theta_0\geq 0$. Then for every $T>0$ there exists a Young measure $\nu$ and functions $(u,\theta)$ defined on $(0,T)$ satisfying the following 
	equalities:
	$$
	\int_0^T\int_a^bu_{tt}v+\int_0^T\int_a^bu_xv_x-\mu\int_0^T\int_a^b\theta_x v=0
	$$
	and
	$$
	\int_0^T\int_a^b\theta_t\psi+\int_0^T\int_a^b\theta_x\psi_x+\mu\int_0^T\int_a^bu_t\theta\psi_x=-\mu\int_0^T\int_a^b\langle \nu_{t,x};w\gamma\rangle\psi ,
	$$
	for all test functions $(v,\psi)\in H_0^1(a,b)\times H^1(a,b)$, and
	$$
	\langle\nu_{t,x};w\rangle:=\int_{(a,b)^2}wd\nu_{t,x}(w) =u_t (t,x) ,\quad
	\langle\nu_{t,x};\gamma\rangle:=\int_{(a,b)^2}\gamma d\nu_{t,x}(\gamma)=u_x (t,x).
	$$
	Moreover, the solution $(u,\theta)$ satisfies the first order estimates:
	\begin{equation}
	\|u\|_{W^{1,\infty}(0,T;L^2(a,b))}+\|u\|_{L^{\infty}(0,T;H^1(a,b))}+\|\theta\|_{L^{\infty}(0,T;L^1(a,b))}\leq C\Big (\|u_0\|_{H^1(a,b)},\|u_1\|_{L^2(a,b)},\|\theta_0\|_{L^2(a,b)}\Big )
	\end{equation}
	and
	\begin{equation}
	\|\theta\|_{L^{\infty}(0,T;L^2(a,b))}+\|\theta\|_{L^2(0,T;H^1(a,b))}\leq C\Big (\|u_0\|_{H^1(a,b)},\|u_1\|_{L^2(a,b)},\|\theta_0\|_{L^2(a,b)}\Big ).
	\end{equation}
	
\end{theorem}
\fi

\begin{theorem}\label{thm:exist_global}
	Let $u_0\in H_0^1(a,b)$, $u_1\in L^2(a,b)$, $\theta_0\in L^2(a,b)$, $\theta_0\geq 0$. 
	Then for every $T>0$ there exists a measure $\gamma\in L^2(0,T;{\mathcal M}(a,b))$
	and functions $(u,\theta)$ defined on $(0,T)$ satisfying the following 
	equalities:
	$$
	\int_0^T\int_a^bu_{tt}v+\int_0^T\int_a^bu_xv_x+\mu\int_0^T\int_a^b\theta_x v=0
	$$
	and
	$$
	\int_0^T\int_a^b\theta_t\psi+\int_0^T\int_a^b\theta_x\psi_x-\mu\int_0^T\int_a^b
	\theta u_t\psi_x=\mu\int_0^T\int_a^b\psi d\gamma,
	$$
	for all test functions $(v,\psi)\in H_0^1(a,b)\times H^1(a,b),$
	with $(u,\theta)$ and $\gamma$ being related in the following way: 
	there exists a sequence $(u^n,\theta^n)$ such that
	\begin{align*}
	(u_t^n,\theta_x^n)\rightharpoonup (u_t,\theta_x)\;{\rm weakly\; in}\; L^2(0,T;L^2(a,b)),\\
	u_t^n\theta_x^n\rightharpoonup \gamma\;{\rm weakly\; in}\; L^2(0,T;\mathcal{M}(a,b)).
	\end{align*}
	Moreover, there exist a constant $C>0,$ depending on the initial data, 
	such that the solution $(u,\theta)$ satisfies the following first order estimates:
	$$
	\|u\|_{W^{1,\infty}(0,T;L^2(a,b))}+\|u\|_{L^{\infty}(0,T;H^1(a,b))}+\|\theta\|_{L^{\infty}(0,T;L^1(a,b))}
	\leq C\Big (\|u_0\|_{H^1(a,b)},\|u_1\|_{L^2(a,b)},\|\theta_0\|_{L^2(a,b)}\Big )
	$$
	and
	$$
	\|\theta\|_{L^{\infty}(0,T;L^2(a,b))}+\|\theta\|_{L^2(0,T;H^1(a,b))}
	\leq C\Big (\|u_0\|_{H^1(a,b)},\|u_1\|_{L^2(a,b)},\|\theta_0\|_{L^2(a,b)}\Big ).
	$$
	
\end{theorem}

\section{Literature review}

It is well known that the equations of one-dimensional nonlinear thermoelasticity
in general admit smooth, classical solution; locally for any data and globally
for small data. This was investigated by various authors and one of the pioneer
works was Slemrod's paper \cite{slemrod} from 1981. 

He considered the following nonlinear thermoelasticity problem 
\begin{align}\label{slem}
\begin{split}
u_{tt}-a(u_x,\theta)u_{xx}+b(u_x,\theta)\theta_x=0,\\
c(u_x,\theta)\theta_t - d(u_x,\theta)\theta_{xx} + b(u_x,\theta)u_{tx}=0,
\end{split}
\end{align}
with mixed boundary conditions (Neumann for $u$, Dirichlet for $\theta$ and vice versa),
where $a,b,c,d$ are differentiable functions such that $a,c,d>0$, $|b|>0$. 
The author proved local existence and uniqueness of a solution.
Furthermore, assuming the smallness of the initial data,
he also proved global existence and uniqueness.
We emphasize here that a strong assumption that $b$ is bounded away from $0$
prevents an expected non-negativity result for $\theta$.
The methodology of the proof is based on application of 
the contraction mapping theorem to solutions of 
a related linear problem. 

Later on, Racke \cite{racke} proved the local existence of classical smooth solutions to the
equations of one-dimensional nonlinear thermoelasticity \eqref{slem}
for the physically reasonable Dirichlet boundary conditions
(the boundary of the configuration is assumed to be rigidly clamped 
and held at constant temperature) for 
both bounded and unbounded domains assuming smooth data.

Racke and Shibata dealt with the same problem in \cite{rashi}
where they proved a global existence of 
smooth solutions by using the spectral analysis to estimate the
decay rates of solutions to the linearized problem, which are then
used in a standard way to obtain energy estimates. Racke and Shibata
\cite{rashi} as well as Racke \cite{racke} again assume $|b|>0$.
Consequently, temperature $\theta$ constructed by them does not share
any comparison principle and the non-negativity result for temperature is missing.
Our model is the simplest one addressing the potential non-negativity of temperature.
Indeed, our solutions are proven to be non-negative.

By following the similar approach as in \cite{rashi}, the global existence 
of smooth solutions for small and smooth data in the case of
Neumann boundary conditions was proved in \cite{Shibata}.

The results of \cite{rashi} were further improved by Racke, Shibata and Zheng in \cite{rsz}. 
They prove that if the initial data are close to the equilibrium then the problem admits
a unique, global, smooth solution. This is achieved without assuming additionally that 
$b$ is bounded away from $0$. However, this was done in the case of Dirichlet boundary data. 
Moreover, the required regularity of initial data ($u_0,\theta_0\in H^3$) 
is much higher than the one we need. Last, but not least, no temperature non-negativity results are 
presented in \cite{rsz}. However, since the construction of the solution there does not require
$b$ to be bounded away from $0$, one could hope for the comparison principle 
in the parabolic equation which would yield non-negativity of temperature. This would clearly require 
some additional work. Not fully trivial, due to the presence of the zeroth order term 
in the parabolic equation, one would have to make sure that the obtained regularity 
of $u_{tx}$ is high enough.

Let us also mention \cite{jiangracke}, where the authors prove a local existence 
theorem for quasilinear system of thermoelasticity in 3D with Dirichlet boundary conditions. 
Again, $b$ is not assumed to be bounded away from $0$, but the regularity requirement concerning 
initial data is much higher than in our case. Again, one could address the positivity of the temperature 
in that case via comparison principle for the parabolic equation. Similarly as in \cite{rsz}, it
is not straightforward. Handling the presence of the zeroth order term would again need quite high 
regularity of $u_{tx}$, hence additional work.

Some more general models were also considered. For instance, global existence of solutions for small initial data and decay of
classical solutions for the equations of one-dimensional nonlinear thermoelasticity
was also considered by Hrusa and Tarabek \cite{HrusaTara},  Jiang \cite{jiang1990} and
Zheng and Shen \cite{Zheng}.

In their paper \cite{HuWang}, Hu and Wang investigated the global solvability of
smooth small solutions to the one-dimensional thermoelasticity problem with second sound in 
the half line. Their work was motivated by Jiang's paper \cite{jiang1993} where the author 
obtained the global existence of smooth solutions to the system of classical thermoelasticity
under Dirichlet boundary conditions in the half line by directly using energy estimates.

Regarding singularities, 
Dafermos and Hsiao in \cite{DafHsiao} considered a special one-dimensional 
model taking into account the whole real line as a reference configuration.
They proved that for large data a smooth solution blows up in a finite time.
Similar problem was also considered by Hrusa and Messaoudi \cite{HrusaMess}.
They showed the existence of smooth initial data for which the solution will
develop singularities in finite time.

Comparing our work to the existing literature, we see that we assume
less regularity on the initial data which leads to a different functional framework.
Moreover, we do not assume that $b(u_x,\theta)$ is bounded away from zero,
which is the case in e.g. \cite{racke,rashi,slemrod}.
Furthermore, none of the works mentioned above dealt with the non-negativity of
the temperature, which is, besides the proof of existence and uniqueness
of the solution, one of the main novelties of our paper.

\section{Preliminaries}\label{sec:prel}

In the paper, we will repeatedly use the following interpolation inequality 
(often called \emph{Agmon's inequality}), stated and proven below 
(in one-dimensional framework) for completeness.
\begin{proposition}\label{prop:agmon}
For any $f\in H_0^1(a,b)$ the following inequality holds 
\begin{equation}\label{agmon_zero}
\|f\|_{L^{\infty}}\leq \sqrt{\|f\|_{L^2}\|f_x\|_{L^2}},
\end{equation}
while if only $f\in H^1(a,b)$, then
\begin{equation}\label{agmon_raz}
\|f\|_{L^{\infty}}\leq \sqrt{2\|f\|_{L^2}\|f_x\|_{L^2}}+\sqrt{\frac{1}{b-a}\int_a^b f^2(s)ds}.
\end{equation}
\end{proposition}
\proof
First we notice that for any $f\in C^1[a,b]$ and $x_0<x$ one has
\begin{equation}\label{raz}
f^2(x)=\int_{x_0}^x (f^2)'(s) ds +f^2(x_0),
\end{equation}
hence for $x_0$ such that the value of $f^2$ at $x_0$ is smaller than the average
\[
f^2(x)\leq 2\int_{x_0}^xf(s)f'(s)ds+f^2(x_0)\leq 2\|f\|_{L^2}\|f_x\|_{L^2}+\frac{1}{b-a}\int_a^b f^2(s)ds
\]
where we used the Cauchy-Schwarz inequality. Next, the usual density argument allows us to 
state \eqref{agmon_raz} for any $f\in H^1(a,b)$. In order to arrive at \eqref{agmon_zero} 
we pick up first $x_0$ being $a$ in \eqref{raz} to arrive at
\[
f^2(x)=\int_a^x \left(f^2(s)\right)' ds,
\]
next for $x_0$ being $b$ in \eqref{raz} we have
\[
f^2(x)=\int_b^x \left(f^2(s)\right)' ds.
\] 
Hence, for smooth functions $f$ supported in $(a,b)$ one has
\[
2f^2(x)\leq 2\left(\int_a^x |f(s)||f'(s)|ds + \int_x^b |f(s)||f'(s)|ds\right)=2\int_a^b |f(s)||f'(s)|ds.
\]
By applying the Cauchy-Schwarz inequality to the right-hand side we obtain
\[
f^2(x)\leq \|f\|_{L^2}\|f_x\|_{L^2}, 
\]
and using the density of compactly supported smooth functions in $H_0^1(a,b)$ we arrive at \eqref{agmon_zero}.

\qed

Moreover, in Section~\ref{sec:aux} we shall need a version of the Schaefer's theorem
for the subspace of nonnegative functions.
Since we could not find a proper reference we attach a version of this theorem 
(together with the proof) which is applicable in our case.
\begin{proposition}\label{prop:schaef}
Let $\bar{X}$ be the Banach space of real-valued functions, by $X$ we denote its subset $\bar{X}\cap\{f\geq 0 \}$. 
Let $A:X\rightarrow X$ be compact and continuous mapping, assume moreover that the set of all points for which there exists $\lambda \in [0,1]$ so that 
\begin {equation}\label{dwa}
B:=\{x\in X: x=\lambda Ax\}
\end{equation}
is bounded. Assume next that the topology in $\bar{X}$ is order preserving 
(closure of the set of non-negative functions consists of non-negative functions). Then $A$ has a fixed point $x\in X$.
\end{proposition}
\proof
The proof follows the lines of the original proof of Schaefer (see \cite{schaefer} or compare to \cite{evans}). 
We only have to ensure that the fixed point is a non-negative function. After choosing a constant $M$ in such a way that 
\[
\left\|u\right\|< M\;\;\mbox{for all}\;\; u\in B\;(\mbox{the definition of B is given in} \eqref{dwa}),
\]
one can define $A_M(u):=\frac{MA(u)}{\left\|A(u)\right\|}$ for $u$ such that $\left\|A(u)\right\|\geq M$, while 
$A_M=A$ for $u$ such that $\left\|A(u)\right\|\leq M$. We observe that 
$A_M: B(0,M)\cap X\rightarrow B(0,M)\cap X$,
$A_M$ inherits continuity and compactness from $A$. Hence 
\[
A_M:\overline{conv (B(0,M)\cap X)}\rightarrow \overline{conv (B(0,M)\cap X)}.
\] 
We are in a position to use classical Schauder's fixed point theorem and say that 
$A_M$ has a fixed point in $\overline{conv (B(0,M)\cap X)}$.
On the one hand the latter set consists of non-negative functions only. 
Indeed, a convexification of the set of non-negative functions consists of non-negative functions only.
Non-negativity is also preserved by the closure, as assumed in the statement of the proposition.
Finally, we show that a fixed point of $A_M$ is also a fixed point of $A$ in a standard way.    

\qed

\section{Auxiliary problem}\label{sec:aux}

In this section we introduce an auxiliary problem which we shall utilize in order to construct solutions of
\eqref{Evolution}. More precisely, we follow the \emph{artificial viscosity approach} (see e.g. \cite{ChenDaf}), i.e. 
we find a global solution to the regularized system:
\begin{equation}\label{Viscoelastic}
\begin{array}{c}
u_{tt}-u_{xx}-\nu u_{txx}=-\mu \theta_x\quad {\rm in}\quad (0,T)\times (a,b), \\
\theta_t-\theta_{xx}+\mu\theta(u_t)_x=0,\;\theta\geq 0
\quad{\rm in}\quad (0,T)\times (a,b),\\
u(t,a)=u(t,b)=0,\; \theta_x(t,a)=\theta_x(t,b)=0,\\
u(0,.)=u_0,\; u_t(0,.)=u_1,\\
\theta(0,.)=\theta_0\geq 0,
\end{array}
\end{equation}
and  then construct a solution to the system \eqref{Evolution}
as a limit of solutions to \eqref{Viscoelastic} when $\nu\rightarrow 0,$ $\nu>0$ being a regularization parameter (artificial structural viscoelasticity). The initial conditions
$\theta_0, u_0$ and $u_1$ belong to the following function spaces:
\begin{equation}\label{IC}
u_0\in H_0^1(a,b), \; u_1\in L^2(a,b), \; \theta_0\in L^2(a,b).
\end{equation}

\begin{remark}
Full thermoviscoelastic system has also an additional term in the temperature equation
(see e.g. \cite{Dafsima}). Equation \eqref{Viscoelastic}$_2$ with physical viscosity would be:
$$
\theta_t-\theta_{xx}+\mu\theta(u_t)_x=\nu (u_{tx})^2.
$$
\end{remark}

\begin{definition}\label{def:visco_sol}
We say that $(u,\theta)$ is a weak solution to problem \eqref{Viscoelastic} if the following conditions are satisfied:
\begin{enumerate}
\item 
\[
u\in W^{1,\infty}(0,T;L^2(a,b))\cap H^1(0,T; H_0^1(a,b)),\quad 
u_{tt}\in L^2(0,T;H^{-1}(a,b))
\]
\item  
\[ 
\theta \in H^1(0,T; L^2(a,b))\cap L^2(0,T;H^2(a,b)),\quad
\theta\geq 0.
\]
\item $(u,\theta)$ satisfies the following variational equation:
\begin{equation}\label{weak1}
_{H^{-1}}\langle u_{tt},v\rangle_{H^1} +\int_a^bu_xv_x+\nu\int_a^bu_{tx}v_x+\mu\int_a^b\theta_x v=0,
\end{equation}
\begin{equation}\label{weak2}
\int_a^b\theta_t\psi+\int_a^b\theta_x\psi_x+\mu\int_a^b\theta u_{tx}\psi=0,
\end{equation}
where $(v,\psi)\in H_0^1(a,b)\times H^1(a,b).$
\item Initial conditions \eqref{Viscoelastic}$_{5,6}$ are satisfied.
\end{enumerate}
\end{definition}

\iffalse
\begin{definition}\label{def:visco_sol}
Function $u$ with the following regularity
\[
u\in W^{1,\infty}(0,T;L^2(a,b))\cap H^1(0,T; H_0^1(a,b))
\]
satisfying moreover
\[
u_{tt}, u_{xx}\in L^2(0,T;L^2(a,b))
\]
and a nonnegative function 
\[
\theta \in H^1(0,T; L^2(a,b))\cap L^2(0,T;H^2(a,b))
\]
are weak solutions to problem \eqref{Viscoelastic} if
\begin{equation}\label{weak1}
\int_a^bu_{tt}v dx+\int_a^bu_xv_xdx+\nu\int_a^bu_{tx}v_xdx-\mu\int_a^b\theta_x vdx=0
\end{equation}
\begin{equation}\label{weak2}
\int_a^b\theta_t\psi+\int_a^b\theta_x\psi_x-\mu\int_a^b\theta u_{tx}\psi=0,
\end{equation}
where $(v,\psi)\in H_0^1(a,b)\times H^1(a,b).$
\end{definition}
\fi
\begin{remark}
Due to required regularity of the temperature $\theta\in L^2(0,T;H^1(a,b))$ 
as well as $\theta_t\in L^2(0,T;L^2(a,b)),$ 
we have that $$\theta\in C([0,T];L^2(a,b))$$ (see e.g. \cite{evans}). 
Similarly, from the regularity of the displacement $u,$ we get that:
\begin{equation*}
u\in C([0,T]; H^1(a,b)) \text{ and } u_t\in C([0,T];L^2(a,b)).
\end{equation*}
Therefore solutions belong to a regularity class in which initial conditions \eqref{Viscoelastic}$_{5,6}$ can be understood in the strong sense.
\end{remark}

We will prove the existence of global-in-time solution $(u,\theta)$ of \eqref{Viscoelastic} 
by using the second order energy estimate and the Schaefer's fixed point theorem as 
introduced in Section~\ref{sec:prel}. Let us choose an arbitrary (but fixed) $T>0$ and 
first define suitable function space:
\begin{equation}\label{SpaceFPT1}
\mathcal{H}(0,T)=\{\theta\in L^2(0,T;H^1_0(a,b)):\theta\geq 0\}.
\end{equation}
The strong topology in $\mathcal{H}$ satisfies the order-preserving assumption in Proposition~\ref{prop:schaef}.

Next we define operator $F:\mathcal{H}(0,T)\to \mathcal{H}(0,T) $ in the following way. 
Take $\tilde{\theta}\in \mathcal{H}(0,T)$ and define  $\tilde{u}$ as
a solution of the following initial-boundary value problem:
\begin{equation}\label{WaveTilde}
\left\{
\begin{array}{l}
\tilde{u}_{tt}-\tilde{u}_{xx}-\nu \tilde{u}_{txx}=-\mu \tilde{\theta}_x\; {\rm in}\;(0,T)\times (a,b),\\
\tilde{u}(.,a)=\tilde{u}(.,b)=0\; {\rm on}\; (0,T),\\
\tilde{u}(0,.)=u_0,\; \tilde{u}_t(0,.)=u_1\; {\rm on}\; (a,b).
\end{array}
\right.
\end{equation}
Now, with given $\tilde{u}$, we define $\theta$ as a solution of the following
initial-boundary value problem:
\begin{equation}\label{ThetaTilde}
\left\{
\begin{array}{l}
\theta_t-\theta_{xx}=-\mu\theta\tilde{u}_{tx}\; {\rm in}\; (0,T)\times (a,b),\\
\theta_x(.,a)=\theta_x(.,b)=0 \; {\rm on}\; (0,T),\\
\theta(0,.)=\theta_0\geq 0 \;{\rm on}\;  (a,b).
\end{array}
\right.
\end{equation}
We define $F(\tilde{\theta}):=\theta$.

\begin{theorem}\label{thm:first}
Let $u_0\in H_0^1(a,b), u_1\in L^2(a,b),\theta_0\in H^1(a,b),\theta_0\geq 0$. 
Then $F:\mathcal{H}(0,T)\to \mathcal{H}(0,T)$ is continuous. Moreover, we have the 
following estimates for $\theta=F(\tilde{\theta})$
\begin{equation}\label{RegEstimates}
\|\theta\|_{H^1(0,T;L^2(a,b))}+\|\theta\|_{L^{2}(0,T;H^2(a,b))}\leq C\|\tilde{\theta}\|_{\mathcal{H}(0,T)},
\end{equation}
where $C$ depends only on initial data and $T$.
\end{theorem}
\proof
We will prove the theorem in series of lemmas.
\begin{lemma}\label{lem:wave_est}
There exists a unique solution to problem \eqref{WaveTilde} with the following properties:
\begin{equation*}
\|\tilde{u}\|_{L^{\infty}(0,T;H^1(a,b))}+ \|\tilde{u}\|_{W^{1,\infty}(0,T;L^2(a,b))}+\nu\|\tilde{u}\|_{H^1(0,T;H^1(a,b))}\leq C\|\tilde{\theta}\|_{\mathcal{H}(0,T)}.
\end{equation*}
\end{lemma}
\proof
Since \eqref{WaveTilde} has a fixed right-hand side, it is a linear equation. Moreover, it is a linear damped wave equation (for $\tilde{u}$). The construction of a solution can be done in a standard way by using the Galerkin approximations (see e.g. \cite{evans}), and therefore here we just derive the  formal energy estimate. We take $\tilde{u}_t$ as a test function to obtain:
\begin{align*}
\frac{1}{2}\frac{d}{dt}(\|\tilde{u}_t\|^2_{L^2}+\|\tilde{u}_x\|^2_{L^2})+\nu\|\tilde{u}_{tx}\|^2_{L^2}\leq\frac{1}{2}(\|\tilde{u}_t\|^2_{L^2}+\|\tilde{\theta}_x\|^2_{L^2}).
\end{align*}
Now the uniqueness and the required estimates follow directly from Gronwall's inequality.

\iffalse
In order to arrive at the required regularity for $\tilde{u}_{tt}$ in Definition~\ref{def:visco_sol} we take $\tilde{u}_{tt}$ as a test function to obtain:
\begin{align*}
\frac{\nu}{2}\|\tilde{u}_{tx}(t)\|^2_{L^2}
+\|\tilde{u}_{tt}\|^2_{L^2_tL^2_x}
&=\|\tilde{u}_{tx}\|^2_{L^2_tL^2_x}-\int_a^bu_x(t)u_{tx}(t)+\int_a^b u_0 (u_1)_x
+\frac{\nu}{2}\|(u_1)_x\|^2_{L^2}
+\nu\int_0^t\int_a^b\tilde{\theta}_x \tilde{u}_{tt}
\\
&\leq \|\tilde{u}_{tx}\|^2_{L^2_tL^2_x}+\|u_x(t)\|_{L^2}\|u_{xt}(t)\|_{L^2}
+\|u_0\|_{H^1}\|u_1\|_{H^1}
+\frac{\nu}{2}\|(u_1)_x\|^2_{L^2}
+\|\tilde{\theta}_x\|_{L^2_tL^2_x}\|\tilde{u}_{tt}\|_{L^2_tL^2_x}.
\end{align*}
\fi

\qed

\begin{remark}
	Throughout the rest of the manuscript we will use $\|\cdot\|_{L^2}$ to denote
	$L^2$-norm in space (and analogously for any other norm).
	For the norms in space and time, the shortened notation
	$\|\cdot\|_{L^2_t L^2_x}$ will be used (and analogously for any other norm).
\end{remark}

\begin{lemma}\label{lem:heat_est}
There exists a unique solution to problem \eqref{ThetaTilde} with the following properties:
\begin{equation*}
\|\theta\|_{H^1(0,T;L^2(a,b))}+\|\theta\|_{L^{2}(0,T;H^2(a,b))}\leq C \|\tilde{\theta}\|_{\mathcal{H}(0,T)}.
\end{equation*}
\end{lemma}
\proof
This lemma follows the classical results in parabolic equations (see \cite{lsu}, Chapter III).
 More precisely, since $\tilde{u}_{tx}\in L^2_t(L^2_x)$, the existence and uniqueness of a
 weak solution is a consequence of Theorem~3.1 (p. 145) and Theorem~4.1 (p. 153). Moreover, 
 Theorem~7.1 (p. 181) and Corollary~7.1 (p. 186) imply that $\theta\in L^{\infty}_t(L^{\infty}_x)$ 
 and therefore $\theta u_{tx}\in L^2_t(L^2_x)$. Now, the statement of the lemma follows from Theorem~6.1 (p. 178).

\qed

\begin{lemma}[Positivity]\label{lem:posit}
$\theta\geq 0$ on $(0,T)\times (a,b)$.
\end{lemma}
\proof
Let us consider the following problem
$$
\theta_t-\theta_{xx}=-\mu\theta_+\tilde{u}_{tx},
$$
where $\theta_+=\max\{\theta,0\}$ is a positive part of $\theta$. By taking $\theta_-=\max\{-\theta,0\}$ as a test function we obtain:
$$
-\frac{1}{2}\frac{d}{dt}\|\theta_-\|^2_{L^2}-\|(\theta_-)_x\|^2_{L^2}=0.
$$
Since $\theta_-(0,.)=0$, we obtain that $\theta_-=0$. Positivity of $\theta$ now follows from Lemma~\ref{lem:heat_est} which guarantees the uniqueness of the solution of \eqref{ThetaTilde}.

\qed

\begin{lemma}[Continuity]\label{lem:cont}
$F:\mathcal{H}(0,T)\to \mathcal{H}(0,T)$ is continuous.
\end{lemma}
\proof
The continuity follows from Lemma~\ref{lem:wave_est} and Lemma~\ref{lem:heat_est}.
Let us take sequence $\tilde{\theta}_n\to \tilde{\theta}$ in $\mathcal{H}(0,T)$ and denote by $\tilde{u}_n$, $\tilde{u}$ the solutions of \eqref{WaveTilde} corresponding to $\tilde{\theta}_n$, $\tilde{\theta}.$ Set $\theta_n=F(\tilde{\theta}_n)$, $\theta=F(\tilde{\theta})$.

Since problem \eqref{WaveTilde} is linear, functions $v_n=\tilde{u}_n-\tilde{u}$ satisfy the following equation with zero initial and boundary conditions:
$$
(v_n)_{tt}-(v_n)_{xx}-\nu (v_n)_{txx}=-\mu (\tilde{\theta}_n-\tilde{\theta})_x\quad {\rm in}\; (0,T)\times (a,b).
$$
By taking $(v_n)_t$ as a test function, similarly as in Lemma~\ref{lem:wave_est}, we conclude $(\tilde{u}_n)_{tx}\to \tilde{u}_{tx}$ in $L^2(0,T;L^2(a,b))$. Now, let $\psi_n=\theta_n-\theta$. Then functions $\psi_n$ satisfy the following equation with zero initial and boundary data:
$$
(\psi_n)_t-(\psi_n)_{xx}+\mu\psi_n\tilde{u}_{tx}=-\mu\theta((\tilde{u}_n)_{tx} - \tilde{u}_{tx}).
$$
From Lemma~\ref{lem:heat_est} it follows that $\theta$ is bounded and therefore the right-hand side of the last equation converges to zero in $L^2(0,T;L^2(a,b))$. The statement of the lemma follows from the stability of a weak solution (see e.g. \cite{lsu}, Theorem~4.5 (p. 166)).

\qed \\
By summing up all the results obtained in the previous lemmas, we see that 
Theorem~\ref{thm:first} is proved. 
\begin{theorem}\label{thm:visco_exist}
Let $u_0\in H_0^1(a,b), u_1\in L^2(a,b),\theta_0\in H^1(a,b)$, $\theta_0\geq 0$. 
Then for any $T>0$ there exists a solution $(u,\theta)$ to problem \eqref{Viscoelastic} 
on $(0,T)$ satisfying the following estimate:
\begin{equation}\label{ViscSolEst}
\|u\|_{W^{1,\infty}(0,T;L^2(a,b))}+\|u\|_{H^1(0,T;H^1(a,b))}+
\|\theta\|_{H^1(0,T;L^2(a,b))}+\|\theta\|_{L^{2}(0,T;H^2(a,b))}
\leq C(\nu,{\rm data}).
\end{equation}
\end{theorem}
\proof
We use Schaefer's fixed point theorem (Proposition~\ref{prop:schaef}). 
From Theorem~\ref{thm:first} we conclude that $F$ is a continuous and compact mapping. 
This follows from the fact that the set
$$
\{f\in \mathcal{H}(0,T): \|f\|_{\mathcal{H}(0,T)}+\|f\|_{H^1(0,T;L^2(a,b))}+\|f\|_{L^2(0,T;H^2(a,b))}\leq C\}
$$
is compact in $\mathcal{H}(0,T)$ as a consequence of Aubin-Lions lemma for the triplet 
$$
H^2(a,b)\subset\subset H^1(a,b)\subset L^2(a,b).
$$ 
It remains to prove the boundedness of the following set:
$$
\{f\in \mathcal{H}(0,T):f=\lambda F(f),\; \lambda\in [0,1]\}.
$$
Let $\lambda\in [0,1]$ and $\theta=\lambda F(\theta)$. Then $\theta$ is a solution to system \eqref{Viscoelastic}. Because of Lemma~\ref{lem:heat_est}, equation is satisfied in strong sense and therefore we can integrate equation \eqref{Viscoelastic}$_2$ over $(a,b).$ Furthermore, we take $u_t$ as a test function in \eqref{Viscoelastic}$_1$ and sum the resulting equations to obtain the following estimate:
$$
\frac{d}{dt}\left (\frac{1}{2}\|u_t\|^2_{L^2}+\frac{1}{2}\|u_x\|^2_{L^2}+\|\theta\|_{L^1}\right )+\nu\|u_{tx}\|^2_{L^2}=0.
$$
Therefore, we conclude that $\|u_{tx}\|_{L^2(0,T;L^2(a,b))}\leq C$, where $C$ depends only on the initial data. Now, the bound for $\theta$ can be obtained by using the standard results for parabolic equation \cite{lsu} analogously as in Lemma~\ref{lem:heat_est}.

\qed

\begin{proposition}\label{prop:unique1}
	Weak solution of the regularized problem \eqref{Viscoelastic} obtained in Theorem \ref{thm:visco_exist} is unique.
\end{proposition}
\proof
We start a proof by recalling that due to \eqref{ViscSolEst}, there exists a constant $C>0$, which depends on $\nu$ and data, such that 
\begin{equation}\label{tak}
\left\|\theta_x\right\|_{L^2(0,T;H^1(a,b))}<C, 
\end{equation}
and
\[
\left\|\theta_t\right\|_{L^2(0,T;L^2(a,b))}<C.
\]
From the last inequality we infer
\begin{equation}\label{taak}
\left\|\theta_{tx}\right\|_{L^2(0,T;H^{-1}(a,b))}\leq C\left\|\theta_t\right\|_{L^2(0,T;L^2(a,b))}<C.
\end{equation}
In view of \eqref{tak} and \eqref{taak}, we arrive at (see \cite{evans})
\begin{equation}\label{taaak}
\theta_x\in C([0,T];L^2(a,b)).
\end{equation} 

Let $(u_1,\theta_1)$, $(u_2,\theta_2)$ be two weak solutions of problem 
\eqref{Viscoelastic} and  set $u=u_1-u_2$, $\theta=\theta_1-\theta_2$. 
By subtracting \eqref{Viscoelastic}$_1$ for $(u_1,\theta_1)$ and $(u_2,\theta_2)$ 
we get that $u$ satisfies the following differential equation with zero initial and boundary conditions:
$$
u_{tt}-u_{xx}-\nu u_{txx}=-\mu \theta_x.
$$
By multiplying the above equation by $u_t$, integrating over space and time interval, and using Young's and Gronwall's inequalities, we get:
$$
\|u_t\|_{L_t^{\infty}L_x^2}+\|u_x\|_{L_t^{\infty}L_x^2}+\nu\|u_{tx}\|_{L_t^{2}L_x^2}\leq C\|\theta_x\|_{L^2_tL^2_x}.
$$
The equation for $\theta=\theta_1-\theta_2$ reads:
$$
\theta_t-\theta_{xx}=-\mu(\theta_1u_{tx}+\theta (u_2)_{tx}).
$$
We multiply the above equation by $\theta_t$ and integrate over $(a,b)$ to obtain:
\begin{align*}
\|\theta_t\|^2_{L^2}+\frac{1}{2}\frac{d}{dt}\|\theta_x\|^2_{L^2}
&=-\mu\int_a^b(\theta_1u_{tx}\theta_t+\theta(u_2)_{tx}\theta_t)\\
&\leq C\left (\|\theta_1\|_{L^{\infty}}\|u_{tx}\|_{L^2}\|\theta_t\|_{L^2}
+\|\theta\|_{L^{\infty}}\|(u_2)_{tx}\|_{L^2}\|\theta_t\|_{L^2}\right)
\end{align*}
for a.e. $t\in(0,T).$
From the proof of Lemma~\ref{lem:heat_est} we know that $\|\theta_1\|_{L^{\infty}_tL^{\infty}_x}\leq C.$ 
Having that in mind, 
we use Young's inequality to bound the right-hand side of the previous inequality:
$$
\|\theta_t\|^2_{L^2}+\frac{1}{2}\frac{d}{dt}\|\theta_x\|^2_{L^2}
\leq C\varepsilon\|\theta_t\|^2_{L^2}+\frac{C}{\varepsilon}\|u_{tx}\|^2_{L^2}
+\frac{C}{\varepsilon}\|\theta\|^2_{L^2}\|(u_2)_{tx}\|_{L^2}^2
+\frac{C}{\varepsilon}\|\theta_x\|^2_{L^2}\|(u_2)_{tx}\|_{L^2}^2,
$$
where we also used that $\|\theta\|_{L^\infty}^2\leq c\left(\|\theta\|_{L^2}^2+\|\theta_x\|_{L^2}^2\right)$.
By choosing $\varepsilon$ such that $C\varepsilon< \frac{1}{2},$
the first term on the right-hand side can be absorbed into 
the first term on the left-hand side, and  after integrating from $0$ to $t$ we obtain:
\begin{align*}
\frac{1}{2}\left(\int_0^t\|\theta_t(s)\|^2_{L^2}ds+\|\theta_x(t)\|_{L^2}^2\right)
\leq 
\frac{C}{\varepsilon}\|\theta\|^2_{L_t^\infty L^2_x}\int_0^t\|(u_2)_{tx}\|_{L^2}^2+\frac{C}{\varepsilon}\|u_{tx}\|^2_{L^2_tL^2_x} + \frac{C}{\varepsilon}\int_0^t\|\theta_x\|^2_{L^2}\|(u_2)_{tx}\|_{L^2}^2.
\end{align*}
Multiplying the inequality with 2 and using the estimate $\|u_{tx}\|_{L_t^2L_x^2}\leq C\|\theta_x\|_{L_t^2L_x^2}$ we obtain:

\begin{equation}\label{Uniq1}
\int_0^t\|\theta_t(s)\|^2_{L^2}ds+\|\theta_x(t)\|_{L^2}^2\leq 
\frac{C}{\varepsilon}\|\theta_{x}\|^2_{L^2_tL^2_x} +\frac{C}{\varepsilon}\|\theta(t)\|^2_{L^2}\int_0^t\|(u_2)_{tx}\|_{L^2}^2+ \frac{C}{\varepsilon}\int_0^t\|\theta_x\|^2_{L^2}\|(u_2)_{tx}\|_{L^2}^2.
\end{equation}

Next, using $\theta(0)=0$, we have
\[
\left\|\theta(t)\right\|_{L^2}=\left\|\int_0^t\theta_t(s)ds \right\|_{L^2}\leq \sqrt{t}\sqrt{\int_0^t\left\|\theta_t\right\|_{L^2}^2}.
\]
Hence, \eqref{Uniq1} turns into
\begin{align}\label{Uniq2}
\begin{split}
&\int_0^t\|\theta_t(s)\|^2_{L^2} ds+\|\theta_x(t)\|_{L^2}^2\\
&\leq \frac{C}{\varepsilon}\|\theta_{x}\|^2_{L^2_tL^2_x} +\frac{C}{\varepsilon}t\int_0^t\left\|\theta_t\right\|_{L^2}^2\int_0^t\|(u_2)_{tx}\|_{L^2}^2+ \frac{C}{\varepsilon}\int_0^t\|\theta_x\|^2_{L^2}\|(u_2)_{tx}\|_{L^2}^2.
\end{split} 
\end{align}

We are now in a position to finish the proof by a variant of an argument 
used in the proof of Lemma~2.1 in \cite[p. 140]{lsu}. 
First we partition interval $(0,T)$ into a finite number of subintervals $(t_{k-1},t_k)$, $k=1,\dots, N$, such that 
\begin{equation}\label{czasy}
\frac{1}{4}\leq\frac{C}{\varepsilon}\int_{t_{k-1}}^{t_k}\|(u_2)_{tx}\|_{L^2}^2\leq \frac{1}{2},\quad t_k-t_{k-1}<1,\quad k=1,\dots,N.
\end{equation}
This can be done because $\|u_{tx}\|_{L^2}$ is square integrable.   %$\|\theta_x(t^*)\|_{L^2}=\|\theta_x\|_{L^{\infty}L^2}>0$. 

Now we proceed inductively, first we prove that $\theta_x=0$ on $(t_0=0,t_1)$. Keeping in mind \eqref{czasy}, from \eqref{Uniq2} we have
\begin{equation}\label{main}
\frac{1}{2}\int_{0}^{t}\|\theta_t(s)\|^2_{L^2}ds+\|\theta_x(t)\|_{L^2}^2
\leq \frac{C}{\varepsilon}\int_0^{t}\|\theta_{x}\|^2_{L^2}
+\frac{C}{\varepsilon}\|\theta_x\|^2_{L^{\infty}_tL^2_x}\int_0^{t}\|(u_2)_{tx}\|_{L^2}^2
\end{equation}
for any $t\leq t_1$.

Let us fix $0<t<t_1$. Due to \eqref{taaak}, we notice that 
$f(t):=\underset{0\leq s\leq t}{\sup}\left\|\theta_x(s,x)\right\|_{L^2(a,b)}$ 
is bounded, in particular integrable in $(0,t)$.
In view of \eqref{main}, for $0<s<t$ we have
\begin{align*}
\|\theta_x(s)\|_{L^2}^2&\leq 
\frac{C}{\varepsilon}\int_0^{s}\|\theta_{x}\|^2_{L^2}
+\frac{C}{\varepsilon}\sup_{0\leq z\leq s}\|\theta_x(z,x)\|^2_{L^2}\int_0^s\|(u_2)_{tx}\|_{L^2}^2\\
&\leq \frac{C}{\varepsilon}\int_0^t\sup_{0\leq z\leq s}\|\theta_{x}(z,x)\|^2_{L^2}ds
+\frac{C}{\varepsilon}\sup_{0\leq s\leq t}\|\theta_x(s,x)\|^2_{L^2}\int_0^t\|(u_2)_{tx}\|_{L^2}^2.
\end{align*} 
Next, we take supremum of the left-hand side over $s\leq t$ 
\[
\sup_{0\leq s\leq t}\|\theta_x(s)\|_{L^2}^2\leq \frac{C}{\varepsilon}
\int_0^t\sup_{0\leq z\leq s}\|\theta_{x}(z,x)\|^2_{L^2}ds
+\frac{C}{\varepsilon}\sup_{0\leq s\leq t}\|\theta_x(s,x)\|^2_{L^2}\int_0^t\|(u_2)_{tx}\|_{L^2}^2.
\]
Finally, utilizing \eqref{czasy}, the last inequality turns into
\[
\frac{1}{2}\sup_{0\leq s\leq t}\|\theta_x(s)\|_{L^2}^2
\leq \frac{C}{\varepsilon}\int_0^t\sup_{0\leq z\leq s}\|\theta_{x}(z,x)\|^2_{L^2}ds.
\]
Using Gronwall's inequality we have
\begin{equation}\label{waz}
\theta_x=0\quad\rm{on}\quad [0,t_1].
\end{equation} 
Next,  Poincar\'e's inequality tells us that $\theta(t,x)=c(t)$. 
On the other hand, \eqref{main} shows, in view of \eqref{waz}, 
that $\theta=0$ on $[0,t_1]$. In particular, $\theta(t_1)=0$ and 
therefore we can further iterate the argument to prove that $\theta=0$ 
on $[t_{k-1},t_k]$, $k=1,\dots,N$, and thus finish the uniqueness proof.

\qed

We shall next 
prove that our unique solution is actually regular provided initial data is more regular.
Indeed, if we impose more restrictive assumptions on initial data, Theorem~\ref{thm:visco_exist}
still gives a unique solution. We shall give formal estimates which can be consequently used to
arrive at more regular solutions via a Schaefer's theorem like in the proof of Theorem \ref{thm:visco_exist}.
This time, since we require more regularity from our solutions we search for a fixed point in the set
\[
{\cal H}_1(0,T)=\{\theta\in H^1(0,T;H^1_0(a,b):\theta\geq 0\}.
\] 
In such a case an obtained solution overlaps with the solution constructed in Theorem \ref{thm:visco_exist}, see the uniqueness claim in Proposition \ref{prop:unique1}. The obtained solution is smoother, see a proposition below.  
\begin{proposition}\label{prop:visco_reg}
Let $u_0\in H^2(a,b)\cap H^1_0(a,b)$, $u_1\in H^2(a,b)\cap H^1_0(a,b)$, 
$\theta_0\in H^3(a,b)$, $\theta_0\geq 0$. Then the unique solution $(u,\theta)$ 
to problem \eqref{Viscoelastic}
given by Theorem~\ref{thm:visco_exist} satisfies the following regularity properties:
\begin{equation}\label{ViscSolReg}
\|u\|_{W^{2,\infty}(0,T;L^2(a,b))}+\|u\|_{H^2(0,T;H^1(a,b))}+
\|\theta\|_{H^2(0,T;L^2(a,b))}+\|\theta\|_{L^{2}(0,T;H^3(a,b))}
\leq C(\nu,{\rm data}).
\end{equation}

\end{proposition}
\proof
First of all, by Theorem \ref{thm:visco_exist}, for any $t>0$ we have
\begin{equation}\label{utx}
\int_0^t\int_a^b u_{tx}^2 dx< C.
\end{equation}
Next, multiplying \eqref{Viscoelastic}$_2$ by $\theta_t$ and then integrating in space and time,
using \eqref{utx}, we arrive at
\begin{equation}\label{thetat1}
\int_0^t\int_a^b \theta_t^2 dxdt <C\;\mbox{for any}\;\; t>0. 
\end{equation}   
Next, we apply $\partial_t$ to \eqref{Viscoelastic}$_1$, multiply the resulting equation by $u_{tt}$ 
and integrate in both space and time to get:
$$
\|u_{tt}\|_{L^{\infty}_tL^2_x}+\|u_{tx}\|_{L^{\infty}_tL^2_x}+2\nu\|u_{ttx}\|_{L^2_tL^2_x}
\leq C\|\theta_{tx}\|_{L^2_tL^2_x}+\|u_{tt}(0)\|_{L^2}+\|u_{tx}(0)\|_{L^2}.
$$
Since, $u_{tt}=u_{xx}+\nu u_{txx}-\mu\theta_x$, we conclude:
\begin{equation}\label{ViscoHighEst}
\|u_{tt}\|_{L^{\infty}_tL^2_x}+\|u_{tx}\|_{L^{\infty}_tL^2_x}+2\nu\|u_{ttx}\|_{L^2_tL^2_x}
\leq C\big (\|\theta_{tx}\|_{L^2_tL^2_x}+\|u_0\|_{H^2}+\|u_1\|_{H^2}+\|\theta_0\|_{H^1}\big ).
\end{equation}
Now, we apply $\partial_t$ to both sides, multiply the resulting equation by $\theta_{tt}$ and integrate to get:
$$
\|\theta_{tt}\|_{L^2_tL^2_x}^2+\frac{1}{2}\|\theta_{tx}(t)\|_{L^2}^2
=-\mu\int_0^t\int_a^b \theta_t u_{tx}\theta_{tt}-\mu\int_0^t\int_a^b\theta u_{ttx}\theta_{tt}
+\frac{1}{2}\|\theta_{tx}(0)\|_{L^2}^2=\rm{I+II+III}.
$$
Let us estimate terms on the right-hand side separately. 
Using \eqref{agmon_raz} and Young's inequality we get:
\begin{align*}
|\rm{I}|&\leq \int_0^t\|u_{tx}\|_{L^{2}}\|\theta_t\|_{L^{\infty}}\|\theta_{tt}\|_{L^2}
\leq \int_0^t\|u_{tx}\|_{L^{2}}(\|\theta_t\|^{1/2}_{L^2}\|\theta_{tx}\|^{1/2}_{L^2}
+\|\theta_t\|_{L^2})\|\theta_{tt}\|_{L^2}\\
&\leq C\varepsilon\|\theta_{tt}\|_{L^2_tL^2_x}^2+\frac{C}{\varepsilon}
\int_0^t\|u_{tx}\|_{L^2}^2(\|\theta_t\|^2_{L^2}+\|\theta_{tx}\|^2_{L^2})\\
&\leq C\varepsilon\|\theta_{tt}\|_{L^2_tL^2_x}^2+\frac{C}{\varepsilon}\left(
\sup_{0\leq s\leq t}\|u_{tx}(s)\|_{L^2}^2\int_0^t\|\theta_t\|^2_{L^2}+\int_0^t\|u_{tx}\|_{L^2}^2\|\theta_{tx}\|^2_{L^2}\right)\\
&\leq C\varepsilon\|\theta_{tt}\|_{L^2_tL^2_x}^2+
\frac{C}{\varepsilon}
\left(\int_0^t \left\|u_{tx}\right\|_{L^2}^2\|\theta_{tx}\|_{L^2}^2+\left(\int_0^t\|\theta_{tx}\|_{L^2}^2+\|u_0\|_{H^2}+\|u_1\|_{H^2}+\|\theta_0\|_{H^1}\right)\int_0^t\|\theta_t\|^2_{L^2}\right),
\end{align*}
where in the last inequality we used \eqref{ViscoHighEst}.
Consequently, in view of \eqref{utx} and \eqref{thetat1}
$$
|\text{I}|\leq C\varepsilon\|\theta_{tt}\|_{L^2_tL^2_x}^2+\frac{C}{\varepsilon}
\int_0^t \left\|u_{tx}\right\|_{L_2}^2\|\theta_{tx}\|_{L^2}^2+ C(\|u_0\|_{H^2},\|u_1\|_{H^2},\|\theta_0\|_{H^1})\int_0^t\|\theta_{tx}\|_{L^2}^2.
$$
Moreover, from the proof of Lemma~\ref{lem:heat_est} 
we have $\|\theta\|_{L^{\infty}_tL^{\infty}_x}\leq C$. Therefore using \eqref{ViscoHighEst} we have:
\begin{align*}
|\rm{II}|&\leq C\varepsilon\|\theta_{tt}\|_{L^2_tL^2_x}^2+\frac{C}{\varepsilon}\|u_{ttx}\|_{L^2_tL^2_x}^2\\ 
&\leq C \varepsilon
\|\theta_{tt}\|_{L^2_tL^2_x}^2+\frac{C}{\varepsilon}\left(\|\theta_{tx}\|_{L^2_tL^2_x}^2+\|u_0\|_{H^2}+\|u_1\|_{H^2}+\|\theta_0\|_{H^1}\right).
\end{align*}
Finally, 
$$
|\text{III}|\leq C(\|\theta_{xxx}(0)\|_{L^2}^2+\|(\theta_x u_{tx})(0)\|_{L^2}^2+\|(\theta u_{txx})(0)\|_{L^2}^2)\leq C(\nu,{\rm data}).
$$
Picking up $\varepsilon$ small enough and summing up the estimates of I, II and III, we arrive at
\begin{equation}\label{koncowe}
\|\theta_{tx}(t)\|_{L^2}^2\leq \frac{C}{\varepsilon}
\int_0^t \left\|u_{tx}\right\|_{L_2}^2\|\theta_{tx}\|_{L^2}^2+C(\varepsilon, {\rm data})\int_0^t\|\theta_{tx}\|_{L^2}^2+C(\nu,\varepsilon, {\rm data}).
\end{equation}
We are now in a position to finish the proof by using slightly less standard version of Gronwall's inequality.
For reader's convenience, we give the details. First, we simplify the notation and rewrite \eqref{koncowe} as
\[
\|\theta_{tx}(t)\|_{L^2}^2\leq C(\varepsilon, {\rm data})
\int_0^t (\left\|u_{tx}(s)\right\|_{L_2}^2+1)\|\theta_{tx}(s)\|_{L^2}^2ds+C(\nu,\varepsilon, {\rm data}).
\]
Taking into account \eqref{utx}, multiplying both sides by $(\|u_{tx}(t)\|_{L^2}^2+1)e^{-\int_0^t (\|u_{tx}(s)\|_{L^2}^2+1)ds}$, the latter can be transformed into
\[
\frac{d}{dt}\left(\int_0^t \|\theta_{tx}(s)\|_{L^2}^2(\|u_{tx}(s)\|_{L^2}^2+1)ds\ e^{-\int_0^t(\|u_{tx}(s)\|_{L^2}^2+1)ds}\right)\leq
-C(\nu,\varepsilon, {\rm data})\frac{d}{dt}e^{-\int_0^t(\|u_{tx}(s)\|_{L^2}^2+1)ds},
\]
and consequently after integration
\[
\int_0^t \|\theta_{tx}(s)\|_{L^2}^2(\|u_{tx}(s)\|_{L^2}^2+1)ds\leq C(\nu,\varepsilon, {\rm data})\left(e^{\int_0^t\|u_{tx}(s)\|_{L^2}^2ds}e^t-1\right)
\]
and hence, in view of \eqref{koncowe}
\begin{align*}
\|\theta_{tx}(t)\|_{L^2}^2 &\leq C(\varepsilon)
\int_0^t \left(\left\|u_{tx}\right\|_{L_2}^2+1\right)\|\theta_{tx}\|_{L^2}^2 +C(\nu,\varepsilon, {\rm data})\\
&\leq
C(\nu,\varepsilon, {\rm data})\left(e^{\int_0^t\|u_{tx}(s)\|_{L^2}^2ds}e^t-1\right)+C(\nu,\varepsilon, {\rm data}).
\end{align*}
\qed

\section{Time-independent estimates}\label{sec:TIest}

In this section we derive the time-independent estimates for the solution
$(u,\theta)$ of \eqref{Viscoelastic}. We will use them to construct 
global-in-time weak solutions to \eqref{Evolution}. 
 
We first get the same estimate as in the proof of Theorem~\ref{thm:visco_exist}. 
We multiply \eqref{Viscoelastic}$_1$ by $u_t$ and integrate from $a$ to $b$, 
next integrate \eqref{Viscoelastic}$_2$, and sum the resulting expressions 
to obtain the basic energy equality:
\begin{equation}\label{BasicEnergy}
\frac{d}{dt}\left(\frac{1}{2}\|u_t\|^2_{L^2}+
	\frac{1}{2}\|u_x\|^2_{L^2}+\|\theta\|_{L^1}\right)+\nu \|u_{tx}\|^2_{L^2}=0
\end{equation}
which yields the following estimate:
\begin{equation}\label{BasicEst1}
\|u\|_{W^{1,\infty}(0,T;L^2(a,b))}+\|u\|_{L^{\infty}(0,T;H^1(a,b))}+\|\theta\|_{L^{\infty}(0,T;L^1(a,b))}\leq C\left(\|u_t(0)\|_{L^2},\|u_x(0)\|_{L^2},\|\theta(0)\|_{L^2}\right).
\end{equation}

To derive the second estimate we multiply \eqref{Viscoelastic}$_2$ by $\theta$ and integrate:
$$
\frac{1}{2}\frac{d}{dt}\|\theta\|^2_{L^2}+\|\theta_x\|^2_{L^2}
=-\mu\int_a^b\theta^2(u_t)_x =2\mu\int_a^b\theta_x\theta u_t 
\leq 2\mu\|\theta\|_{L^{\infty}}\|\theta_x\|_{L^2}\|u_t\|_{L^2}.
$$
Now, from \eqref{BasicEst1} we conclude that $\|u_t(t)\|_{L^2}$ is bounded for a.e. $t\in [0,T]$. 
Next, \eqref{agmon_raz} and Young's inequality yield
\begin{align*}
\frac{1}{2}\frac{d}{dt}\|\theta\|^2_{L^2}+\|\theta_x\|^2_{L^2}
&\leq C_1\|\theta\|^{1/2}_{L^2}\|\theta_x\|_{L^2}^{3/2}
+ C_2\|\theta\|_{L^2}\|\theta_x\|_{L^2}\\
&\leq \left(\frac{C_1}{\varepsilon^3}\|\theta\|_{L^2}^2+C_1\varepsilon\|\theta_x\|^2_{L^2} \right)
+ \left(\frac{C_2}{\varepsilon}\|\theta\|_{L^2}^2+C_2\varepsilon\|\theta_x\|^2_{L^2} \right)\\
&\leq  \frac{C}{\varepsilon^3}\|\theta\|_{L^2}^2+C\varepsilon\|\theta_x\|^2_{L^2},
\end{align*}
where $C=2\mu\max\{\sqrt{2},\sqrt{\frac{1}{b-a}}\}.$
We choose $\varepsilon$ such that the last term can be absorbed into the right-hand side and use Gronwall's inequality to obtain:
\begin{equation}\label{BasicEst2}
\|\theta\|_{L^{\infty}(0,T;L^2(a,b))}+\|\theta\|_{L^2(0,T;H^1(a,b))}
\leq \exp(Ct)\|\theta(0)\|_{L^2}.
\end{equation}

\section{Time-dependent estimates}\label{sec:TDest}

This section is devoted to higher order estimates of solutions to \eqref{Viscoelastic}. 
Those estimates will hold only on properly short time intervals. 
We will need them to obtain local-in-time well-posedness of \eqref{Evolution}.
The following theorem holds true:
\begin{theorem}\label{thm:short_time} 
	Let $u_0\in H^2(a,b)\cap H^1_0(a,b)$, $u_1\in H^2(a,b)\cap H^1_0(a,b)$, $\theta_0\in H^3(a,b)$. 
	Then there exists short enough time $T_0>0$ and a constant $C>0$ such that $(u,\theta)$, 
	solutions of \eqref{Viscoelastic} given by Proposition~\ref{prop:visco_reg}, satisfy
\begin{align}\label{gorne}
\begin{split}
&\|u\|_{W^{2,\infty}(0,T_0;L^2(a,b))}+\|u\|_{W^{1,\infty}(0,T_0;H^1(a,b))}
+\|u\|_{L^{\infty}(0,T_0;H^2(a,b))}\\
&\quad\leq C\left(\|u_{tt}(0)\|_{L^2(a,b)}+\|u_{tx}(0)\|_{L^2(a,b)}+\|\theta_t(0)\|_{L^2(a,b)}\right)
\end{split}
\end{align}
and
\begin{equation}\label{dolne}
\|\theta\|_{W^{1,\infty}(0,T_0;L^2(a,b))}+\|\theta\|_{H^1(0,T_0;H^1(a,b))}
\leq C\left(\|u_{tt}(0)\|_{L^2(a,b)}+\|u_{tx}(0)\|_{L^2(a,b)}+\|\theta_t(0)\|_{L^2(a,b)}\right).
\end{equation}
\end{theorem}

\proof
The proof of the theorem will consist of a few steps.

\noindent
\textbf{Step 1.}
We apply $\partial_t$ to \eqref{Viscoelastic}$_1$, next multiply the outcome by $u_{tt}$ 
and integrate over $(a,b)$ to obtain:
$$
\frac{1}{2}\frac{d}{dt}\|u_{tt}\|^2_{L^2}-\int_a^bu_{txx}u_{tt}-\nu\int_a^bu_{txxt}u_{tt}
\leq \frac{1}{2}\left(\|u_{tt}\|^2_{L^2}+\mu\|\theta_{tx}\|^2_{L^2}\right).
$$
Next, integrating by parts the second and third term on the left-hand side 
and utilizing the boundary conditions $u_{tt}(a)=u_{tt}(b)=0,$ being a consequence
 of the fact that the value of $u$ at those points is equal to zero 
(by \eqref{Viscoelastic}$_3$), we arrive at
$$
\frac{1}{2}\frac{d}{dt}\left(\|u_{tt}\|^2_{L^2}+\left\|u_{tx}\right\|^2_{L^2}\right)+\nu\left\|u_{ttx}\right\|^2
\leq \frac{1}{2}\left(\|u_{tt}\|^2_{L^2}+\mu\|\theta_{tx}\|^2_{L^2}\right).
$$
Dropping the $\nu$-term and using Gronwall's inequality we arrive at 
\begin{equation}\label{Energy2Tmp1}
\|u_{tt}\|_{L^{\infty}_t L^2_x}+\|u_{tx}\|_{L^{\infty}_t L^2_x}
\leq  \exp(t)\left(\|u_{tt}(0)\|_{L^2}+\|u_{tx}(0)\|_{L^2}+\mu\|\theta_{tx}\|_{L^2_t L^2_x}\right).
\end{equation}

\noindent
\textbf{Step 2.} Let us now differentiate \eqref{Viscoelastic}$_2$ with respect to time, 
multiply the resulting equation by $\theta_t$ and integrate over $(a,b)$ to obtain:
\begin{equation*}
\frac{1}{2}\frac{d}{dt}\|\theta_t\|^2_{L^2}+\|\theta_{tx}\|^2_{L^2}
=-\mu\int_a^b \left(\theta_t^2(u_t)_x + \theta\theta_t(u_{tt})_{x}\right).
\end{equation*}
The right-hand side can be rewritten by using integration by parts and 
H\"{o}lder's inequality:
\begin{align*}
\frac{1}{2}\frac{d}{dt}\|\theta_t\|^2_{L^2}+\|\theta_{tx}\|^2_{L^2}
&=-\mu\int_a^b(\theta_t^2(u_t)_x-u_{tt}(\theta_x\theta_t+\theta\theta_{tx}))\\
&\leq \mu\|u_{tx}\|_{L^2}\|\theta_t\|_{L^2}\|\theta_t\|_{L^{\infty}}
+\mu\|u_{tt}\|_{L^2}\|\theta_x\|_{L^{\infty}}\|\theta_t\|_{L^2}
+\mu\|u_{tt}\|_{L^2}\|\theta\|_{L^{\infty}}\|\theta_{tx}\|_{L^2}.
\end{align*}
Notice that the boundary terms coming from integration by parts vanish due to
$u_{tt}(a)=u_{tt}(b)=0.$
We integrate the resulting inequality from $0$ to $t,$  
use \eqref{Energy2Tmp1} to estimate $\|u_{tx}\|_{L^2}$ and $\|u_{tt}\|_{L^2}$
and use Agmon's inequality \eqref{agmon_zero}, \eqref{agmon_raz} to obtain:
\begin{align}\label{thetat}
\begin{split}
&\frac{1}{2}\|\theta_t(t)\|_{L^2}^2+\|\theta_{tx}\|^2_{L^2_t L^2_x}\\ 
&\leq \frac{1}{2}\|\theta_t(0)\|_{L^2}^2+\mu\exp(t)
\left(\|u_{tt}(0)\|_{L^2}+\|u_{tx}(0)\|_{L^2}+\mu\|\theta_{tx}\|_{L^2_t L^2_x}\right)\\
&\quad\int_{0}^t\left(\|\theta_t\|_{L^{\infty}}\|\theta_t\|_{L^2}+\|\theta_x\|_{L^{\infty}}\|\theta_t\|_{L^2} 
+\|\theta\|_{L^{\infty}}\|\theta_{tx}\|_{L^2}\right) \\
&\leq \frac{1}{2}\|\theta_t(0)\|_{L^2}^2+c\mu\exp(t)
\left(\|u_{tt}(0)\|_{L^2}+\|u_{tx}(0)\|_{L^2}+\mu\|\theta_{tx}\|_{L^2_t L^2_x}\right)\\
&\quad \int_{0}^t \left(\|\theta_t\|_{L^2}^{3/2}\|\theta_{tx}\|_{L^2}^{1/2}+\|\theta_t\|_{L^2}^2
+\|\theta_x\|^{1/2}_{L^2}\|\theta_{xx}\|_{L^2}^{1/2}\|\theta_t\|_{L^2}
+\|\theta\|_{L^2}^{1/2}\|\theta_x\|_{L^2}^{1/2}\|\theta_{tx}\|_{L^2}
+\|\theta\|_{L^2}\|\theta_{tx}\|_{L^2}\right).
\end{split}
\end{align}
We set 
\begin{equation}\label{aaa}
A(t):=c\mu^2\exp(t),
\end{equation}
where $c=\max\left\{\sqrt{2},\sqrt{\frac{1}{b-a}}\right\},$ 
and firstly estimate the integral in \eqref{thetat}
multiplied by $A(t)\|\theta_{tx}\|_{L^2_t L^2_x}.$

Using H\"older's and Young's inequality (with $p=4/3$ and $q=4$), 
we estimate the first term:
\begin{align}\label{term1}
\begin{split}
\rm{I}&=A(t)\|\theta_{tx}\|_{L_t^2 L_x^2}\int_{0}^t
\|\theta_t(s)\|_{L^2}^{3/2}\|\theta_{tx}(s)\|_{L^2}^{1/2}ds\\
&=A(t)\left(\int_0^t\|\theta_{tx}(s)\|_{L^2}^2ds\right)^{1/2}\left( \int_{0}^t
\|\theta_t(s)\|_{L^2}^{3/2}\|\theta_{tx}(s)\|_{L^2}^{1/2}ds\right)\\
&\leq A(t) \left(\int_0^t\|\theta_{tx}(s)\|_{L^2}^2ds\right)^{1/2}
\left[ \left(\int_0^t\|\theta_t(s)\|_{L^2}^2 ds \right)^{3/4}
\left(\int_0^t\|\theta_{tx}(s)\|_{L^2}^2 ds \right)^{1/4}\right]\\
&= A(t) \left(\int_0^t\|\theta_{tx}(s)\|_{L^2}^2ds\right)^{3/4} \left(\int_0^t\|\theta_t(s)\|_{L^2}^2 ds \right)^{3/4}\\
&\leq A(t)\varepsilon\int_0^t\|\theta_{tx}(s)\|^2_{L^2}ds +\frac{A(t)}{\varepsilon^3}\left(\int_{0}^t\|\theta_t(s)\|_{L^2}^2 ds\right)^3\\
&\leq A(t)\varepsilon\int_0^t\|\theta_{tx}(s)\|^2_{L^2}ds 
+\frac{A(t)\cdot t^2}{\varepsilon^3}\int_0^t\|\theta_t(s)\|_{L^2}^6 ds
\end{split}
\end{align}
In the last inequality, H\"older with $p=\frac{3}{2}$ 
and $q=3$ was used.

%Hence, picking up $\varepsilon$ small enough (such that $A(t)\varepsilon \leq 1/2$), we estimate the first term by
%\begin{equation}\label{czlon_ess}
%\frac{1}{2}\|\theta_{tx}\|^2_{L^2_t L^2_x}+\frac{B(t)}{\varepsilon^3}\int_{0}^t\|\theta_t(s)\|_{L^2}^6 ds
%\end{equation}
%for any $t>0$.

The second term is estimated in a similar way
\begin{align}\label{term2}
\begin{split}
\rm{II}&=A(t)\|\theta_{tx}\|_{L^2_t L^2_x}\int_0^t \|\theta_t(s)\|^2_{L^2} ds\\
&\leq A(t)\varepsilon\int_0^t\|\theta_{tx}(s)\|^2_{L^2}ds +\frac{A(t)}{\varepsilon}\left(\int_{0}^t\|\theta_t(s)\|_{L^2} ds\right)^2\\
&\leq A(t)\varepsilon\int_0^t\|\theta_{tx}(s)\|^2_{L^2}ds 
+\frac{A(t)\cdot t}{\varepsilon}\int_0^t\|\theta_t(s)\|_{L^2}^2 ds.
\end{split}
\end{align}
%By choosing $\varepsilon$ such that $A(t)\varepsilon\leq 1/2,$
%we estimate the first two terms by:
%\begin{equation}
%\frac{1}{2}\|\theta_{tx}\|^2_{L^2_t L^2_x}
%+\frac{A(t)\cdot t^2}{\varepsilon^3}\int_0^t\|\theta_t(s)\|_{L^2}^2 ds.
%\end{equation}

To estimate the third term, we must first estimate $\|\theta_{xx}\|_{L^2},$ which
can be rewritten using \eqref{Viscoelastic}$_2$ and further
estimated by \eqref{agmon_raz}, as follows:
\begin{equation}\label{thetaxx}
\|\theta_{xx}\|_{L^2}\leq \|\theta_t\|_{L^2}+\|\theta u_{tx}\|_{L^2}\leq \|\theta_t\|_{L^2}+\|\theta\|_{L^2}^{1/2}\|\theta_x\|_{L^2}^{1/2}\|u_{tx}\|_{L^2}+\|\theta\|_{L^2}\|u_{tx}\|_{L^2}.
\end{equation}
Using \eqref{thetaxx} we can rewrite the third term:
\begin{align*}
\begin{split}
\rm{III}
&=A(t)\|\theta_{tx}\|_{L^2_t L^2_x}\int_{0}^t\|\theta_t(s)\|_{L^2}\|\theta_x(s)\|^{1/2}_{L^2}\|\theta_{xx}(s)\|_{L^2}^{1/2}ds \\
&\leq A(t)\|\theta_{tx}\|_{L^2_t L^2_x}\int_{0}^t\left(\|\theta_t(s)\|_{L^2}^{3/2}\|\theta_x(s)\|^{1/2}_{L^2}
+\|\theta(s)\|_{L^2}^{1/4}\|\theta_t(s)\|_{L^2}\|\theta_x(s)\|_{L^2}^{3/4}\|u_{tx}(s)\|^{1/2}_{L^2}\right. \\
&\quad \left.
+\|\theta(s)\|_{L^2}^{1/2}\|\theta_t(s)\|_{L^2}\|\theta_x(s)\|_{L^2}^{1/2}\|u_{tx}(s)\|_{L^2}^{1/2}\right)ds.
\end{split}
\end{align*}
We estimate all the terms from the right-hand side separately by using Young's inequality:
	\begin{equation*}
	\text{III.1} 
	=A(t)\|\theta_{tx}\|_{L^2_t L^2_x}\int_0^t \|\theta_t(s)\|_{L^2}^{3/2}\|\theta_x(s)\|_{L^2}^{1/2} ds
	\leq A(t)\varepsilon\|\theta_{tx}\|_{L^2_t L^2_x}^2 \|\theta_x\|_{L^1_t L^2_x} 
	+ \frac{A(t)}{\varepsilon}\int_0^t \|\theta_t(s)\|_{L^2}^3 ds,
	\end{equation*}
	\begin{align*}
	\rm{III.2}&=A(t)\|\theta_{tx}\|_{L^2_t L^2_x}\int_0^t
	\|\theta(s)\|_{L^2}^{1/4}\|\theta_t(s)\|_{L^2}\|\theta_x(s)\|_{L^2}^{3/4}\|u_{tx}(s)\|^{1/2}_{L^2} ds \\
	&\leq A(t)\|\theta_{tx}\|_{L^2_t L^2_x}\|u_{tx}\|_{L^{\infty}_t L^2_x}^{1/2} 
	\int_0^t \|\theta(s)\|_{L^2}^{1/4}\|\theta_t(s)\|_{L^2}\|\theta_x(s)\|_{L^2}^{3/4} ds\\
	&\leq  A(t) \sqrt{\exp(t)}\|u_{tt}(0)\|_{L^2}^{1/2}\|\theta_{tx}\|_{L^2_t L^2_x} \|\theta\|_{L^{\infty}_t L^2_x}^{1/4} 
	\int_0^t \|\theta_t(s)\|_{L^2}\|\theta_x(s)\|_{L^2}^{3/4} ds\\
	&\quad +
	A(t) \sqrt{\exp(t)}\|u_{tx}(0)\|_{L^2}^{1/2}\|\theta_{tx}\|_{L^2_t L^2_x} \|\theta\|_{L^{\infty}_t L^2_x}^{1/4} 
	\int_0^t \|\theta_t(s)\|_{L^2}\|\theta_x(s)\|_{L^2}^{3/4} ds \\
	&\quad + A^{3/2}(t) \|\theta_{tx}\|_{L^2_t L^2_x}^{3/2}\|\theta\|_{L^{\infty}_t L^2_x}^{1/4} 
	\int_0^t \|\theta_t(s)\|_{L^2}\|\theta_x(s)\|_{L^2}^{3/4} ds\\
	&\leq  A(t)\sqrt{\exp(t)}\|u_{tt}(0)\|_{L^2}^{1/2}
	\varepsilon
	\|\theta_{tx}\|_{L^2_t L^2_x}^{4/3}\|\theta\|_{L^{\infty}_t L^2_x}^{1/3} \|\theta_x\|_{L^1_t L^2_x} 
	+\frac{ A(t)\sqrt{\exp(t)}\|u_{tt}(0)\|_{L^2}^{1/2}}{\varepsilon^3}\int_0^t \|\theta_t(s)\|_{L^2}^4 ds\\
	&\quad + A(t)\sqrt{\exp(t)}\|u_{tx}(0)\|_{L^2}^{1/2}
	\varepsilon
	\|\theta_{tx}\|_{L^2_t L^2_x}^{4/3}\|\theta\|_{L^{\infty}_t L^2_x}^{1/3} \|\theta_x\|_{L^1_t L^2_x} 
	+\frac{ A(t)\sqrt{\exp(t)}\|u_{tx}(0)\|_{L^2}^{1/2}}{\varepsilon^3}\int_0^t \|\theta_t(s)\|_{L^2}^4 ds\\
	&\quad +  A^{3/2}(t)\varepsilon\|\theta_{tx}\|_{L^2_t L^2_x}^{2}\|\theta\|_{L^{\infty}_t L^2_x}^{1/3} \|\theta_x\|_{L^1_t L^2_x} 
	+\frac{A^{3/2}(t)}{\varepsilon^3}\int_0^t \|\theta_t(s)\|_{L^2}^4 ds,
	\end{align*}
	\begin{align*}
	\rm{III.3}&=A(t)\|\theta_{tx}\|_{L^2_t L^2_x}\int_0^t
	\|\theta(s)\|_{L^2}^{1/2}\|\theta_t(s)\|_{L^2}\|\theta_x(s)\|_{L^2}^{1/2}\|u_{tx}(s)\|_{L^2}^{1/2}ds\\
	&\leq  A(t) \sqrt{\exp(t)}\|u_{tt}(0)\|_{L^2}^{1/2}\|\theta_{tx}\|_{L^2_t L^2_x} \|\theta\|_{L^{\infty}_t L^2_x}^{1/2} 
	\int_0^t \|\theta_t(s)\|_{L^2}\|\theta_x(s)\|_{L^2}^{1/2} ds\\
	&\quad +
	A(t) \sqrt{\exp(t)}\|u_{tx}(0)\|_{L^2}^{1/2}\|\theta_{tx}\|_{L^2_t L^2_x} \|\theta\|_{L^{\infty}_t L^2_x}^{1/2} 
	\int_0^t \|\theta_t(s)\|_{L^2}\|\theta_x(s)\|_{L^2}^{1/2} ds \\
	&\quad + A^{3/2}(t) \|\theta_{tx}\|_{L^2_t L^2_x}^{3/2}\|\theta\|_{L^{\infty}_t L^2_x}^{1/2} 
	\int_0^t \|\theta_t(s)\|_{L^2}\|\theta_x(s)\|_{L^2}^{1/2} ds\\
	&\leq  A(t)\sqrt{\exp(t)}\|u_{tt}(0)\|_{L^2}^{1/2}
	\varepsilon
	\|\theta_{tx}\|_{L^2_t L^2_x}^{4/3}\|\theta\|_{L^{\infty}_t L^2_x}^{2/3} \|\theta_x\|_{L^1_t L^2_x}^{2/3}
	+\frac{ A(t)\sqrt{\exp(t)}\|u_{tt}(0)\|_{L^2}^{1/2}}{\varepsilon^3}\int_0^t \|\theta_t(s)\|_{L^2}^4 ds\\
	&\quad + A(t)\sqrt{\exp(t)}\|u_{tx}(0)\|_{L^2}^{1/2}
	\varepsilon
	\|\theta_{tx}\|_{L^2_t L^2_x}^{4/3}\|\theta\|_{L^{\infty}_t L^2_x}^{2/3} \|\theta_x\|_{L^1_t L^2_x}^{2/3}
	+\frac{ A(t)\sqrt{\exp(t)}\|u_{tx}(0)\|_{L^2}^{1/2}}{\varepsilon^3}\int_0^t \|\theta_t(s)\|_{L^2}^4 ds\\
	&\quad + A^{3/2}(t)\varepsilon \|\theta_{tx}\|_{L^2_t L^2_x}^{2}\|\theta\|_{L^{\infty}_t L^2_x}^{2/3} \|\theta_x\|_{L^1_t L^2_x}^{2/3}
	+\frac{A^{3/2}(t)}{\varepsilon^3}\int_0^t \|\theta_t(s)\|_{L^2}^4 ds.
	\end{align*}
Notice that in III.2 and III.3 we made use of \eqref{Energy2Tmp1}. Set 
\begin{equation}\label{bbb}
B(t)=\max\left\{ A(t)\sqrt{\exp(t)}\|u_{tt}(0)\|_{L^2}^{1/2}, A(t)\sqrt{\exp(t)}\|u_{tx}(0)\|_{L^2}^{1/2}, A^{3/2}(t)
\right\}
\end{equation}
and apply \eqref{BasicEst2} to see that for any $t>0$ 
we can estimate the term III by
\begin{equation}\label{term3}
\text{III}\leq B(t)\varepsilon\cdot C^{4/3}(t)\|\theta_{tx}\|^2_{L^2_t L^2_x}
+ \frac{B(t)}{\varepsilon^3}\max\left\{\int_{0}^t\|\theta_t(s)\|_{L^2}^3 ds,
\int_{0}^t\|\theta_t(s)\|_{L^2}^4 ds\right\},
\end{equation}
where 
\begin{equation}\label{ccc}
C(t)=\|\theta(0)\|_{L^2}\exp(2\mu c\,t).
\end{equation}

The fourth and fifth term are estimated by using Cauchy-Schwarz inequality and
\eqref{BasicEst2}:
\begin{align}\label{term4}
\begin{split}
\rm{IV}
&=A(t)\|\theta_{tx}\|_{L^2_tL^2_x}\int_{0}^t\|\theta(s)\|^{1/2}_{L^2}\|\theta_x(s)\|^{1/2}_{L^2}\|\theta_{tx}(s)\|_{L^2} ds \\
&\leq A(t)\|\theta_{tx}\|_{L^2_t L^2_x} \left(\int_0^t \|\theta(s)\|_{L^2}\|\theta_x(s)\|_{L^2} ds\right)^{1/2}
\left(\int_0^t \|\theta_{tx}(s)\|_{L^2}^2 ds\right)^{1/2}\\
&\leq A(t)\|\theta_{tx}\|_{L^2_t L^2_x}^2\|\theta\|_{L^{\infty}_t L^2_x}^{1/2}\left(\int_0^t \|\theta_x (s)\|_{L^2} ds \right)^{1/2}\\
&\leq  A(t)\cdot \sqrt[4]{t}\|\theta_{tx}\|_{L^2_t L^2_x}^2\|\theta\|_{L^{\infty}_t L^2_x}^{1/2} \|\theta_x\|_{L^2_t L^2_x}^{1/2}\\
&\leq A(t)\cdot \sqrt[4]{t} \cdot C(t) \|\theta_{tx}\|_{L^2_t L^2_x}^2,
\end{split}
\end{align}
\begin{align}\label{term5}
\begin{split}
\rm{V}
&=A(t)\|\theta_{tx}\|_{L^2_tL^2_x}\int_{0}^t\|\theta(s)\|_{L^2}\|\theta_{tx}(s)\|_{L^2} ds\\
&\leq A(t)\|\theta_{tx}\|_{L^2_t L^2_x} \left(\int_0^t \|\theta(s)\|_{L^2}^2 ds\right)^{1/2}
\left(\int_0^t \|\theta_{tx}(s)\|_{L^2}^2 ds\right)^{1/2}\\
&\leq A(t)\cdot\sqrt{t}\|\theta_{tx}\|_{L^2_t L^2_x}^2 \|\theta\|_{L^{\infty}_t L^2_x}\\
&\leq A(t)\cdot\sqrt{t} \cdot C(t) \|\theta_{tx}\|_{L^2_t L^2_x}^2.
\end{split}
\end{align}
The terms analogous to I-V having  
$c\mu\exp(t)\|u_{tt}(0)\|_{L^2} $ and $c\mu\exp(t)\|u_{tx}(0)\|_{L^2}$ 
instead of $A(t)\|\theta_{tx}\|_{L^2_t L^2_x}$
are easier to estimate since they contain one power of the critical term 
$\|\theta_{tx}\|_{L^2_t L^2_x}$ less and therefore 
can be directly estimated by using Young's inequality and time-independent estimates.

Finally, by inserting the obtained estimates \eqref{term1}, \eqref{term2}, \eqref{term3}, \eqref{term4}, \eqref{term5} into \eqref{thetat} we arrive at 
\begin{align}\label{Energy2Tmp2}
\begin{split}
&\frac{1}{2}\|\theta_t(t)\|_{L^2}^2+\|\theta_{tx}\|_{L^2_tL^2_x}^2
\leq \frac{1}{2}\|\theta_t(0)\|_{L^2}^2+
D(t) \|\theta_{tx}\|_{L^2_t L^2_x}^2\\ 
&+ E(t) 
\max\left\{\int_0^t \|\theta_t(s)\|_{L^2}^2 ds, 
\int_0^t \|\theta_t(s)\|_{L^2}^3 ds,
\int_0^t \|\theta_t(s)\|_{L^2}^4 ds,
\int_0^t \|\theta_t(s)\|_{L^2}^6 ds\right\}
\end{split}
\end{align} 
with 
$$D(t)=\max\left\{A(t)\varepsilon, B(t)\varepsilon\cdot C^{4/3}(t), A(t)\cdot \sqrt[4]{t} \cdot C(t), A(t) \cdot \sqrt{t}\cdot  C(t) \right\}$$
and
$$E(t):= \max\left\{\frac{A(t) \cdot t^2}{\varepsilon^3}, \frac{A(t)\cdot t}{\varepsilon}, \frac{B(t)}{\varepsilon^3} \right\}, $$
where $A(t)$, $B(t)$ and $C(t)$ are given by \eqref{aaa}, \eqref{bbb}
and \eqref{ccc} respectively.

\noindent
\textbf{Step 3.} In order to prove Theorem~\ref{thm:short_time}, we notice that choosing $\varepsilon$ and short enough time so that $D(t)<1/2$ enables us to absorb the second term from the 
right-hand side in \eqref{Energy2Tmp2} into the left-hand side.
Next, applying Gronwall's inequality yields existence of $T_0$ such that 
the claim of Theorem~\ref{thm:short_time} holds for all $t\in (t_0,T_0)$. 
Notice that the estimate on $\|u\|_{L^{\infty}(0,T_0;H^2(a,b))}$ 
arises as a consequence of the obtained regularity for $u$ and $\theta,$
more precisely, since
$u\in W^{2,\infty}(0,T_0;L^2(a,b))$ and $\theta\in H^1(0,T_0;H^1(a,b)),$
the estimate follows immediately from  \eqref{Evolution}$_1$.

\qed

\iffalse
\noindent
\textbf{Step 4} We can now use the obtained formal estimate to how the regularity of solution of auxiliary problem \eqref{Viscoelastic}. We use proceed analogously as in Section \ref{Aux}. We define function space:
$$
\mathcal{H}_1(0,T)=\{\theta\in H^1(0,T;H^1_0(a,b)):\theta\geq 0\}.
$$
Operator $F$ is defined in the same way as in Section \ref{Aux}. We can differentiate equations \eqref{WaveTilde} and \eqref{ThetaTilde} w.r.t. time variable and use standard for linear wave and heat equation. The formal estimates from steps 1-3 allow to apply
\fi

\section{Local-in-time well-posedness} 
In this section we show that the sequence of solutions to system \eqref{Viscoelastic}, whose existence is guaranteed by Theorem~\ref{thm:visco_exist},
converges to a solution of \eqref{Evolution} as the artificial
viscoelasticity parameter tends to zero.

Assume that the initial data have the following regularity
$u_0\in H^2(a,b)\cap H_0^1(a,b,), u_1\in H_0^1(a,b) $ and $\theta_0\in H^2(a,b),\theta_0\geq 0.$ Set $\nu=\frac{1}{n}$ and let $(u^n,\theta^n)$ be the sequence of solutions to \eqref{Viscoelastic} corresponding to initial data
$(u_0,u_1^n,\theta_0^n)$ introduced by a regularization procedure
described below, where $u_1^n$ and $\theta_0^n$ are the
regularized initial elastic velocity and temperature.

We claim that there exists a sequence $u_1^n$ such that:
\begin{itemize}
	\item[(i)] $u_1^n \in H^2(a,b)\cap H_0^1(a,b),$
	\item[(ii)] $u_1^n \to u_1$ in $H^1(a,b),$
	\item[(iii)] $\nu\|u_1^n\|_{H^2(a,b)}\to 0.$
\end{itemize}
Without loss of generality we firstly assume that the interval $(a,b)$
is symmetric around zero 
(we can always achieve that with an appropriate composition) 
and extend $u_1$ to $\mathbb{R}$ by 0.
Set $c:=\frac{b-a}{2}$ and define
\begin{equation*}
\sigma_n = \frac{c\sqrt{n}+1}{c\sqrt{n}-1}.
\end{equation*}
Then the sequence $\tilde{u}_1^n$ defined in the following way
$$\tilde{u}_1^n(x)=u_1(\sigma_n x).$$	
satisfies $\tilde{u}_1^n=0$ on 
$\mathbb{R}\backslash \left(\frac{a}{\sigma_n},\frac{b}{\sigma_n}\right).$
Finally, set
\begin{equation*}
u_1^n = \tilde{u}_1^n*\eta_n,
\end{equation*}
where $\eta_n$ is the sequence of standard mollifiers (with $\varepsilon=\frac{1}{\sqrt{n}}$)
\begin{equation*}
\eta_n(x) = \sqrt{n}\eta(\sqrt{n}x).
\end{equation*}
%	\begin{equation*}
%	\eta(x)=
%	\left\{
%	\begin{array}{l}
%	\exp \left( \frac{1}{|x|^2-1}\right)\text{ if } |x|<1,\\
%	0 \text{ if } |x|\geq 1.
%	\end{array}
%	\right.
%	\end{equation*}
Since $\frac{1}{\sqrt{n}}<c-\frac{c}{\sigma_n}$
we see that $\text{supp }u_1^n= \left(\frac{a}{\sigma_n},\frac{b}{\sigma_n}\right)
+(-\frac{1}{\sqrt{n}},\frac{1}{\sqrt{n}})$ is a subset of $(a,b).$ Thus $u_1^n\in H_0^1(a,b).$

Statement (ii) follows from the 
definition of mollifiers and the fact that $\underset{n\to\infty}{\lim}\sigma_n=1,$
while statement (iii) is a direct consequence of the fact that $\|u_1^n\|_{H^2(a,b)}$
behaves like $\sqrt{n}.$

We next regularize the initial temperature $\theta_0$ by using  the extension
operator to extend $\theta_0$ to $H^2$-function on the real line 
which we then compose with the mollifiers defined above.

%In Section~\ref{sec:TIest} we have obtained the following uniform estimates
%for $(u^n,\theta^n):$
%\begin{equation*}
%\|u^n\|_{W^{1,\infty}(0,T;L^2(a,b))}+\|u^n\|_{L^{\infty}(0,T;H^1(a,b))}+\|\theta^n\|_{L^{\infty}(0,T;L^1(a,b))}\leq C
%\end{equation*}
%and 
%\begin{equation*}
%\|\theta^n\|_{L^{\infty}(0,T;L^2(a,b))}+\|\theta^n\|_{L^2(0,T;H^1(a,b))}\leq C.
%\end{equation*}
%Using this convergence properties we can pass to the limit in Definition~\ref{def:visco_sol} 
%in all the terms except nonlinear one $\int_a^b\theta^n u^n_{tx}\psi$. 
%In order to pass to the limit in nonlinear term, we will make use of Theorem~\ref{thm:short_time}, which will give us additional
%regularity of the solutions to \eqref{Viscoelastic}.

Finally, we employ Theorem~\ref{thm:short_time} to see that the solutions
of \eqref{Viscoelastic} possess the following regularity for any $t$ smaller than $T_0>0$:
\begin{align}\label{thm5}
\begin{split}
&\|u^n\|_{W^{2,\infty}(0,T_0;L^2(a,b))}+\|u^n\|_{W^{1,\infty}(0,T_0;H^1(a,b))}
+\|u^n\|_{L^{\infty}(0,T_0;H^2(a,b))}+\|\theta^n\|_{W^{1,\infty}(0,T_0;L^2(a,b))}\\
&\quad+\|\theta^n\|_{H^1(0,T_0;H^1(a,b))}
\leq C\left( \|u_{tt}^n(0)\|+\|u_{tx}^n(0)\| +\|\theta_t^n(0)\|_{L^2}\right).
\end{split}
\end{align}
The right-hand side is estimated as follows:
$$
\|u_{tt}^n(0)\|_{L^2}\leq \|u_{xx}^n(0)\|_{L^2}+\nu\|u_{txx}^n(0)\|_{L^2}+\mu\|\theta_x^n(0)\|_{L^2}
\leq C\big (\|u_0\|_{H^2}+\nu\|u_1^n\|_{H^2}+\mu\|\theta_0^n\|_{H^1}\big )\leq C,
$$
$$
\|\theta_t^n(0)\|_{L^2}\leq \|\theta_{xx}^n(0)\|_{L^2}+\|(\theta u_{tx}^n)(0)\|_{L^2}
\leq \|\theta_0^n\|_{H^2} + \| \theta_0^n u_1^n\|_{L^2}
\leq C.
$$
Thus, letting $n\to\infty,$ \eqref{thm5} implies the following convergences:
\begin{align}\label{conv_nonlinear}
\begin{split}
u_{tx}^n\rightharpoonup u_{tx}&\quad {\rm weakly\; in}\quad L^2(0,T_0;L^2(a,b)),\\
\theta^n\to \theta&\quad {\rm in}\quad C([0,T_0];C[a,b]).
\end{split}
\end{align}
We are now in a position to pass to the limit in all the terms of the weak formulation.
Notice that the convergences obtained in \eqref{conv_nonlinear} enable us to pass
to the limit in nonlinear term $\int_a^b \theta^n u_{tx}^n\psi.$
Furthermore,  the boundedness of the sequence $u_{tx}^n$
implies vanishing of the regularization term $\nu \int_a^b u_{tx}^n v_x.$
Therefore we have proved the existence part of 
Theorem~\ref{thm:exist_local} which we state again for completeness:

\begin{thm1}
	Let $u_0\in H^2(a,b)\cap H^1_0(a,b)$, $u_1\in H_0^1(a,b)$, $\theta_0\in H^2(a,b)$, $\theta_0\geq 0$. 
	Then there exists a time $T_0>0$ (depending on data) and a unique solution $(u,\theta)$ 
	to problem \eqref{Evolution} on $(0,T_0)$ with the following regularity:
	\begin{align*}
	&\|u\|_{W^{2,\infty}(0,T_0;L^2(a,b))}+\|u\|_{W^{1,\infty}(0,T_0;H^1(a,b))}
	+\|u\|_{L^{\infty}(0,T_0;H^2(a,b))}\\
	&\quad+\|\theta\|_{W^{1,\infty}(0,T_0;L^2(a,b))}+\|\theta\|_{H^1(0,T_0;H^1(a,b))}\leq C.
	\end{align*}
\end{thm1}

To complete the proof of Theorem1~\ref{thm:exist_local}, we still need to show uniqueness.
The below proposition is devoted to this issue.
\begin{proposition}\label{prop:unique}
	The solution $(u,\theta)$ to problem \eqref{Evolution} given by Theorem~\ref{thm:exist_local} is unique.
\end{proposition}
\proof
Let $(u_1,\theta_1)$, $(u_2,\theta_2)$ be two weak solutions of problem \eqref{Evolution} and set 
$$u=u_1-u_2,\quad\theta=\theta_1-\theta_2.$$
By subtracting \eqref{Evolution}$_1$ for $(u_1,\theta_1)$ and $(u_2,\theta_2)$ we get that $u$ satisfies the following differential equation with zero initial and boundary conditions:
$$
u_{tt}-u_{xx}=-\mu \theta_x\; {\rm in}\; (0,T)\times (a,b).
$$
By multiplying the above equation by $u_t$, integrating over space and time interval, and using Young's and Gronwall's inequality, we get:
\begin{equation}\label{est1}
\|u_t\|_{L_t^{\infty}L_x^2}+\|u_x\|_{L_t^{\infty}L_x^2}
\leq C\|\theta_x\|_{L^2_tL^2_x}.
\end{equation}
The equation for $\theta=\theta_1-\theta_2$ reads:
$$
\theta_t-\theta_{xx}=-\mu(\theta_1u_{tx}+\theta (u_2)_{tx}).
$$
We multiply the above equation by $\theta$ and integrate over $(a,b)$ to obtain:
\begin{align}\label{est2}
\frac{1}{2}\frac{d}{dt}\|\theta\|^2_{L^2}+\|\theta_x\|^2_{L^2}
=-\mu\int_a^b(\theta_1u_{tx}\theta+\theta(u_2)_{tx}\theta)
=-\mu\int_a^b\theta_1u_{tx}\theta-\mu\int_a^b\theta(u_2)_{tx}\theta.
%(\|\theta_1\|_{L^{\infty}}\|u_{tx}\|_{L^2}\|\theta_t\|_{L^2}+\|\theta\|_{L^{\infty}}\|(u_2)_{tx}\|_{L^2}\|\theta_t\|_{L^2}\right)
\end{align}
We integrate the previous equation with respect to time and estimate integrals on the right-hand side separately.
%\begin{align*}
%\left|\int_0^t\int_a^b\theta_1u_{tx}\theta\right|
%&\leq \left|\int_0^t\int_a^b(\theta_1)_x u_t\theta\right| + \left|\int_0^t\int_a^b\theta_1 u_t\theta_x\right|.
%\end{align*}
The first integral is separated into two terms by using integration by parts.
Using \eqref{agmon_raz} we estimate the first term:
\begin{align*}
\left|\int_0^t\int_a^b(\theta_1)_x u_t\theta\right|
&\leq \int_0^t \|(\theta_1)_x\|_{L^2}\|u_t\|_{L^2}\|\theta\|_{L^{\infty}} 
\leq \int_0^t \|(\theta_1)_x\|_{L^2}\|u_t\|_{L^2}\left(\|\theta\|_{L^2}^{1/2}\|\theta_x\|_{L^2}^{1/2}
+C\left\|\theta\right\|_{L^2}\right)\\
&\leq \|(\theta_1)_x\|_{L^{\infty}_t L^2_x} \|u_t\|_{L^{\infty}_tL^2_x}\|\theta\|_{L^2_t L^2_x}^{1/2}
\|\theta_x\|_{L^2_t L^2_x}^{1/2}+C\|(\theta_1)_x\|_{L^{\infty}_t L^2_x} \|u_t\|_{L^{\infty}_tL^2_x}\left\|\theta\right\|_{L_t^2L_x^2}.
\end{align*}
Since $\theta_1$ is a solution, Theorem~\ref{thm:exist_local} implies that
$\|(\theta_1)_x\|_{L^{\infty}_t L^2_x}\leq C.$ Furthermore, from estimate
\eqref{est1} we know that $\|u_t\|_{L_t^{\infty}L_x^2}\leq C \|\theta_x\|_{L^2_t L^2_x},$ 
so we can use Young's inequality (with $p=4$ and $q=4/3$) to see:
\begin{equation*}
\left|\int_0^t\int_a^b(\theta_1)_x u_t\theta\right|
\leq C\varepsilon \|\theta_x\|^2_{L^2_t L^2_x}+
\frac{C}{\varepsilon^3}\|\theta\|^2_{L^2_t L^2_x}.
\end{equation*}

The second term is estimated in a similar way:
\begin{align*}
\left|\int_0^t\int_a^b\theta_1 u_t\theta_x\right|
&\leq \int_0^t \|\theta_1\|_{L^{\infty}}\|u_t\|_{L^2}\|\theta_x\|_{L^2}
\leq \int_0^t  C\varepsilon\|\theta_x\|^2_{L^2} + \frac{C}{\varepsilon}\|u_t\|^2_{L^2}
 \|\theta_1\|^2_{L^\infty}  \\
&\leq C\varepsilon\|\theta_x\|^2_{L^2_t L^2_x} + \frac{C}{\varepsilon}\|u_t\|^2_{L_t^{\infty}L_x^2}
\int_0^t \|\theta_1\|^2_{L^\infty}
\leq C\varepsilon\|\theta_x\|^2_{L^2_t L^2_x} + \frac{C}{\varepsilon}\|\theta_x\|^2_{L^2_t L^2_x}
\int_0^t \|\theta_1\|^2_{L^\infty}.
\end{align*}

What is left is to estimate the second integral in \eqref{est2}:
\begin{align*}
\left|\int_0^t\int_a^b\theta (u_2)_{tx}\theta\right|
&\leq \left|\int_0^t\int_a^b \theta_x (u_2)_t\theta\right| + \left|\int_0^t\int_a^b\theta (u_2)_t\theta_x\right| 
\leq \int_0^t \|\theta\|_{L^2}\|(u_2)_t\|_{L^{\infty}}\|\theta_x\|_{L^2}\\
&\leq \|\theta\|_{L^2_t L^2_x}
\|(u_2)_t\|_{L^{\infty}_tL^\infty_x} \|\theta_x\|_{L^2_t L^2_x}.
\end{align*}
Since $u_2$ is a solution, we have that  $\|(u_2)_t\|_{L^{\infty}_tL^{\infty}_x}\leq C$
so the second integral is estimated by using Young's inequality:
\begin{equation*}
\left|\int_0^t\int_a^b\theta (u_2)_{tx}\theta\right| \leq
C\varepsilon \|\theta_x\|^2_{L^2_t L^2_x} 
+ \frac{C}{\varepsilon} \|\theta\|^2_{L^2_t L^2_x}.
\end{equation*}
At last, we employ the obtained estimates into \eqref{est2} to see that
\begin{equation*}
\|\theta(t)\|_{ L^2}^2 + 2 \|\theta_x\|_{L^2_tL^2_x}^2 \leq 
C\varepsilon \|\theta_x\|^2_{L^2_t L^2_x} + \frac{C}{\varepsilon^3} \|\theta\|^2_{L^2_t L^2_x} + \frac{C}{\varepsilon}\|\theta_x\|^2_{L^2_t L^2_x}
\int_0^t \|\theta_1\|^2_{L^\infty}.
\end{equation*}
We choose $\varepsilon$ small enough that the $\theta_x$ term is absorbed into the left-hand side. Next, similarly as in Proposition~\ref{prop:unique1}, we will use a trick from \cite{lsu}. We partition time interval $(0,T_0)$ into finitely many intervals $(t_{k-1},t_k)$ in such a way that
$\frac{C}{\varepsilon}\int_{t_{k-1}}^{t_k} \|\theta_1\|^2_{L^\infty}<1/2$. Since $\theta_1$ is a solution, \[
\left\|\theta_1\right\|_{L^2(0,T_0;L^\infty(a,b))}<\infty.
\] 
We are thus in a position to proceed inductively as in the proof of Proposition \ref{prop:unique1} and absorb the term 
\[
\frac{C}{\varepsilon}\|\theta_x\|^2_{L^2_t L^2_x}
\int_0^t \|\theta_1\|^2_{L^\infty}
\] 
into the left-hand side at each interval $(t_{k-1}, t_k)$. This allows us
to use Gronwall's inequality at each time interval and obtain $\theta=0$ for $t\in (t_{k-1},t_k)$ for any $k$.
The proof is finished.

\qed

\section{Global-in-time measure valued solution}

In order to prove the existence of a global-in-time weak solution, we would like to pass to the limit as $\nu\to 0,$ i.e. $n\to\infty,$ using only the first order estimates. The only difficulty lies in the nonlinear term:
$$
\mu\int_a^b u^n_{tx}\theta^n\psi=-\mu\int_a^b u^n_t(\theta^n_x\psi+\theta^n\psi_x).
$$
Notice that the above expression is well-defined a.e. in $(0,T)$ because of the following estimate:
$$
\left|\mu\int_a^b u^n_t(\theta^n_x\psi+\theta^n\psi_x)\right|
\leq C \|u_t^n\|_{L^2}\|\theta^n\|_{H^1}\|\psi\|_{H^1}.
$$
Because of the uniform convergence of $\theta^n$ we have $\int_a^b u^n_t\theta^n\psi_x\to \int_a^b u_t\theta\psi_x$. However, from the uniform estimates we can only conclude that $u_t^n\theta_x^n$ is bounded in $L^2_tL^1_x$ and therefore there exists a measure $\gamma$ such that
$$
\int_0^T\int_a^b u^n_t\theta^n_x\psi\to \int_0^T\int_a^b\psi d\gamma.
$$

Therefore we have proved the second main theorem of the manuscript:
\begin{thm2}
	Let $u_0\in H_0^1(a,b)$, $u_1\in L^2(a,b)$, $\theta_0\in L^2(a,b)$, $\theta_0\geq 0$. 
	Then for every $T>0$ there exists a measure $\gamma\in L^2(0,T;{\mathcal M}(a,b))$
	and functions $(u,\theta)$ defined on $(0,T)$ satisfying the following 
	equalities:
	$$
	\int_0^T\int_a^bu_{tt}v+\int_0^T\int_a^bu_xv_x+\mu\int_0^T\int_a^b\theta_x v=0
	$$
	and
	$$
	\int_0^T\int_a^b\theta_t\psi+\int_0^T\int_a^b\theta_x\psi_x-\mu\int_0^T\int_a^b
	\theta u_t\psi_x=\mu\int_0^T\int_a^b\psi d\gamma,
	$$
	for all test functions $(v,\psi)\in H_0^1(a,b)\times H^1(a,b),$
	with $(u,\theta)$ and $\gamma$ being related in the following way: 
	there exists a sequence $(u^n,\theta^n)$ such that
	\begin{align*}
	(u_t^n,\theta_x^n)\rightharpoonup (u_t,\theta_x)\;{\rm weakly\; in}\; L^2(0,T;L^2(a,b)),\\
	u_t^n\theta_x^n\rightharpoonup \gamma\;{\rm weakly\; in}\; L^2(0,T;\mathcal{M}(a,b)).
	\end{align*}
	Moreover, there exist a constant $C>0,$ depending on the initial data, 
	such that the solution $(u,\theta)$ satisfies the following first order estimates:
	$$
	\|u\|_{W^{1,\infty}(0,T;L^2(a,b))}+\|u\|_{L^{\infty}(0,T;H^1(a,b))}+\|\theta\|_{L^{\infty}(0,T;L^1(a,b))}
	\leq C\Big (\|u_0\|_{H^1(a,b)},\|u_1\|_{L^2(a,b)},\|\theta_0\|_{L^2(a,b)}\Big )
	$$
	and
	$$
	\|\theta\|_{L^{\infty}(0,T;L^2(a,b))}+\|\theta\|_{L^2(0,T;H^1(a,b))}
	\leq C\Big (\|u_0\|_{H^1(a,b)},\|u_1\|_{L^2(a,b)},\|\theta_0\|_{L^2(a,b)}\Big ).
	$$
	
\end{thm2}

%\begin{corollary}\label{Young}
%The sequence $u^{\nu_k}_{t}\theta^{\nu_k}_{x}$ generates a Young measure $\nu_{x,k}$
%$$
%\gamma=u_t\theta_x??
%$$
%\end{corollary}
%Indeed, as a consequence of estimate 
%\[
%\int_0^T\int_a^b u^{\nu_k}_{t}\theta^{\nu_k}_{x}\psi dx\leq C
%\]
%and the Markov inequality we see that $u^{\nu_k}_{t}\theta^{\nu_k}_{x}$ is tight,
%so it generates a Young measure. 

\section*{Acknowledgement}
T.C. was partially supported by the OPUS 4 grant number 2012/07/B/ST1/03306. 
T.C. is grateful to K. Che\l{}mi\'nski from Warsaw University of Technology 
for stating the problem and helpful suggestions regarding mechanical foundations
 of the problem as well as the fixed point theorem. 
M.G. and B.M. were partially supported by the Croatian Science Foundation 
(Hrvatska Zaklada za Znanost) grant number IP-2018-01-3706.
This work was partially supported by the grant 346300 for IMPAN 
from the Simons Foundation and the matching 2015-2019 Polish MNiSW fund.

%\newpage

\bibliographystyle{plain}
%\bibliography{thermo_references}	

\end{document}